\documentclass[10 pt]{article}

\usepackage[round]{natbib}
\usepackage[pdftex]{graphicx}
\usepackage{amssymb,amsfonts,amsmath}
\usepackage{enumerate}
\usepackage{setspace}
\usepackage{color}

\newtheorem{Theorem}{Theorem}

\newtheorem{Remark}{Remark}

\newtheorem{Lemma}{Lemma}
\newtheorem{Corollary}{Corollary}

\newcommand{\field}[1]{\mathbf{#1}}
\newcommand{\R}{\field{R}}
\newcommand{\E}{\field{E}}

\newcommand{\N}{\field{N}}
\newcommand{\X}{\field{X}}
\newcommand{\bx}{\field{x}}
\newcommand{\Y}{\field{Y}}
\newcommand{\Lb}{\field{L}}
\newcommand{\I}{\field{I}}
\newcommand{\D}{\field{D}}
\newcommand{\W}{\field{W}}
\newcommand{\Wa}{\field{\bar{W}}}
\newcommand{\Wb}{\field{\tilde{W}}}
\newcommand{\Q}{\field{Q}}
\newcommand{\A}{\field{A}}
\newcommand{\U}{\field{U}}
\newcommand{\betab}{\boldsymbol\beta}
\newcommand{\varepsilonb}{\boldsymbol\varepsilon}

\newcommand{\cp}{\xrightarrow{P}}

\newcommand{\mPVC}[1]{   p(R_{0{#1}}^n|p_{#1},0)   }
\newcommand{\MPVC}[1]{   e^{-p_{#1}a_n}p(R_{{#1}f}^n|M_n,p_{#1})   }

\begin{document}
\date{}
\centerline{\large\bf Linear regression model selection using $p$-values }
\centerline{\large\bf when the model dimension grows}
\bigskip
\centerline{\sc 
By Piotr Pokarowski
\footnote{Institute of Applied Mathematics and Mechanics,
Faculty of Mathematics, Informatics and Mechanics,
University of Warsaw.
E-mail: pokar@mimuw.edu.pl.
},
Jan Mielniczuk
\footnote{Institute of Computer Science, Polish Academy of Sciences and
Faculty of Mathematics and Information Sciences, Warsaw University of Technology.
E-mail: miel@mini.pw.edu.pl.
} and
Pawe\l Teisseyre
\footnote{Institute of Computer Science, Polish Academy of Sciences.
E-mail: teisseyrep@ipipan.waw.pl.}
}
\medskip
Abstract.
We consider  a new criterion-based approach to model selection  in linear regression. Properties  of selection criteria based on $p$-values of a likelihood ratio statistic are studied for  families of linear regression models. We prove that such procedures are consistent 
i.e. the minimal true model is chosen with probability tending to 1 even when the number of models under consideration  slowly increases with a sample size. The simulation study indicates that introduced methods perform promisingly when compared  with  Akaike and Bayesian Information Criteria.
\section{Introduction}
We reconsider a 
problem of  model choice for a linear regression 
\begin{equation}
\label{LM}
\Y = \X\betab + \varepsilonb,
\end{equation}
where $\Y$ is an $n\times 1$ vector of observations which variability we would like to explain, $\X$ is a $n\times M_n$ design matrix  
consisting of vectors of $M_n$
potential regressors collected from $n$ objects 
and $\varepsilonb=(\varepsilon_1,\ldots,\varepsilon_n)'$  is an
unknown vector of  errors, assumed to have $N(0,\sigma^{2}\I)$ distribution.
Vector $\betab=(\beta_1,\ldots,\beta_{M_n})'$ is an unknown vector of parameters. 
In the paper we will consider the  cases  corresponding to experimental and observational data when rows of $\X$ are either deterministic or random.
Suppose that some covariates are 
unrelated to the prediction of $\Y$, so that the corresponding coefficients $\beta_i$ are zero. It is assumed that the true model 
is  a submodel of (\ref{LM}). As it is not a priori known which variables are significant 
in order to make the last assumption realistic it is natural to let the horizon
$M_n$ to grow with $n$ and  allow in this way potentially large models.
\\
Model selection is a core issue of statistical modeling. In a framework of linear regression
the problem  has been intensively studied under various conditions imposed on design matrix ${\bf X}$
and growth of $M_n$. The aim of such procedures is to choose the most parsimonious model describing adequately a given data set.
For the review of these advances we refer to \citet{PoetscherLeeb2008}. 
The main problem  here is a modeler's  dillema  that  underfitting  leads
to omission of important variables in the model whereas overfitting involves unnecessary  parameter estimation for redundant coefficients
which lessens the precision of the model fit.\\
In the article we  contribute to a line of research  in which the  chosen model is the maximiser 
of a chosen criterion function.
In a seminal paper which is typical for this approach  \citet{Akaike1970}, starting with  the idea of maximising the expectation of
 predictive likelihood, has shown that the usual likelihood has to be modified  to obtain an unbiased  estimator of
 the expectation. The likelihood  modified in such a way  is known as Akaike Information Criterion (AIC).
Variety of other modifications of the likelihood followed, with Bayes Information Criterion
 (BIC) being the most frequently used competitor.
Recently, \citet{PokarowskiMielniczuk2010} introduced model selection criteria mPVC and MPVC based on $p$-values of
a likelihood ratio statistic for families of linear models with deterministic covariates and constant dimension. 
The idea in the case of minimal $p$-value criterion mPVC is to consider  the model selection problem from a point of
view of
testing a certain null hypothesis $H_0$ against several hypotheses $H_i$ and to choose the hypothesis (the model)
for which the null hypothesis is most strongly rejected in its favour.
The decision in the case of mPVC is based on  a new criterion which is the minimal $p$-value of the underlying test statistics. 
We stress that the discussed selection method is based on a completely different paradigm than the existing approaches: instead of penalizing the likelihood ratio statistic directly by  subtracting a complexity penalty its appropriate function  is chosen as a selection criterion.\\
We study conditions under which such a rule is consistent i.e. 
 it choses the minimal true model
with probability tending to 1 when the sample size increases.
Our main theoretical result stated in Theorem 1 asserts that  this property holds for the minimal $p$-value criterion mPVC
 provided  $M_n$ increases at a slower rate than $\log n + a_n$ where $a_n$  are weights appearing in the scaling of $p$-values.
 Similar result is proved for maximal $p$-value criterion  MPVC.
Both results apply also to the case when $M_n$ is constant provided the full model (\ref{LM}) is correctly specified. We also introduce and
investigate less computationally demanding greedy versions of the discussed methods. \\
In the last section we present the results of limited simulation study which shows that the introduced methods perform on average better than AIC and BIC criteria. In particular, their performance measured by probability of correct subset detection and prediction error is much more stable when the length of list of models $M_n$ increases i.e. regression model becomes sparse.\\
In the paper we focus mainly on explanation i.e. finding the model which adequately describes the data. Besides the immediate application of
model selection methods to to the second main task of prediction let us mention their use in construction of data-adaptive smooth tests (see e.g. \citet{Ledwina1994}). \\
Problem of linear model selection when the number of possible predictors increases with the sample size has been studied from different angle
by \citet{Shao1997} who defined the optimal submodel to be submodel minimizing the averaged squared prediction error and investigated conditions under which the selected model converges in probability to this model. \citet{Morenoetal2010} considered Bayesian  approach to this problem and proposed using Bayes factors for intrinsic priors as selection criteria.\\
The main contribution of the present paper is establishing consistency of the criteria based on p-values when the linear model dimension grows. The result is proved for the random design as well as for the fixed design scenario, the former being treated in detail. Intrumental in the proofs are Lemmas \ref{Auxilliary}, \ref{LemmaSupmodel}, \ref{LemmaSubmodel} which can be also useful for different purposes.\\
The paper is organized as follows.
In Section 2 we 
introduce considered selection criteria.
 In Section 3  we discuss the imposed assumptions and consistency results for the family of models consisting of
all subsets of  predictors as well as hierarchic family. We also  introduce greedy modifications of the considered criteria.
  Section 4 contains proofs 
of the main results and Section 5 discussion of the results of numerical experiments. Proofs of some auxiliary lemmas are relegated to the Appendix.

\section{Model Selection criteria for linear regression  models based on p-values}
We start by explicitly stating the basic assumption we impose on random-design regression model.
Assume that the rows $\bx_1',\ldots,\bx_n'$ of a  matrix $\X(n\times M_n)$ are iid,
$\bx_l=\bx_l^{(n)}=(x_{l,1}^{(n)},\ldots x_{l,M_n}^{(n)})'$, $l=1,\ldots,n$. 
Throughout we consider the situation  that the minimal true model is fixed i.e. it does not change with $n$.
Vectors  $\{{\bx_1^{(n)}}',\ldots,{\bx_n^{(n)}}'\}$ constitute rows in an array of iid sequences of $M_n$-dimensional
random variables. We impose the condition  that $M_n$ is nondecreasing and that the law of the first $M_n$ coordinates of  $\bx_1^{(n+1)}$ coincides with that of
$\bx_1^{(n)}$ i.e. the  distribution of attributes  considered for a certain sample size remains the same for larger sample sizes.
We also assume throughout that the second moments of coordinates of $\bx_1^{(n)}$ are finite for any $n$.
As any submodel of (\ref{LM}) containing $p_j$ variables $(x_{l,j_1}^{(n)},\ldots, x_{l,j_{p_j}}^{(n)})'$ can be described by set 
of indexes $j=\{j_1,\ldots,j_{p_j}\}$  in order  to make  notation simpler it  will be  referred to as model $j$. 
The minimal true model will be denoted by $t$ and  $p_t$  will be the number of nonzero coefficients in  equation (\ref{LM}). 
The empty model $\Y=\varepsilonb$ will be denoted briefly by $0$ and the full model (\ref{LM}) by $f=\{1,\ldots,M_n\}$.  Note that $M_n=p_f$.
Let $\hat{\betab}_j=(\hat{\beta}_{j_1},\ldots,\hat{\beta}_{j_{p_j}})'$ be a maximum likelihood  (ML) estimator of $\betab$ calculated 
for the considered model $j$. We denote $\hat{\betab}_f$, ML estimator in the full model,  briefly by $\hat{\betab}$.
Let $\cal{M}$ be a certain family of subsets of a set $f$ and
 $\bx_{lt}=(x_{l,t_1}^{(n)},\ldots, x_{l,t_{p_t}}^{(n)})'$ be a vector of variables which pertain to the minimal true model $t$.
Througout this paper with exception of Section 3.2 we will impose the following assumption:\\ (A0) $\E(\bx_{1t}\bx_{1t}')$ is positive definite matrix. \\
The main objective  of model selection is to identify the minimal true model $t$ using data $(\X,\Y)$. Let
${\bf f}_{{\betab},{\sigma}^2}(\Y|\X)$ be the conditional density of $\Y$ given $\X$.
Consider two models $j$ and $k$ where the first model is nested within the second model. Denote by $D_{jk}^{n}$ likelihood ratio test (LRT) 
statistic, based on conditional densities given $\X$, for testing $H_0:$ model $j$ is adequate against hypothesis $H_1:$ model $k$ is adequate whereas $j$ is not, equal to
\begin{equation}
\label{dev}
D_{jk}^{n}=2\log\frac{{\bf f}_{\hat{\betab_k},{\hat\sigma_k}^2}(\Y|\X)}{{\bf f}_{\hat{\betab_j},{\hat \sigma_j}^2}(\Y|\X)},
\end{equation}
where ${\hat\sigma_j}^2=RSS(j)/n$ and $RSS(j)$ is a sum of squared residuals from the  ML fit of the model $j$. We recall that 
ML estimator $\hat{\betab_k}$ coincides with Least Squares estimator of $\betab$.
When $j$ and $k$ are linear models  it turns out that LRT statistic is given explicitly by
\begin{equation*}
D_{jk}^{n}=-n\log\left[\frac{RSS(k)}{RSS(j)}\right]=-n\log(1-R_{jk}^{n}), 
\end{equation*}
where  
\begin{equation}
\label{R}
R_{jk}^{n}=\frac{RSS(j)-RSS(k)}{RSS(j)}
\end{equation}
is coefficient of partial determination of variables belonging to $k\setminus j$ 
given that variables in  set $j$ are  included in the model.
Under the null hypothesis $H_0$ 
it follows from Cochran's theorem (cf. e.g. Section 5.5 in \citet{RencherSchaalje2008}) that  given ${\bf X}$ $RSS(j)\sim \sigma^2\chi^2_j$ and
 $R_{jk}^{n}\sim Beta(\frac{p_k-p_j}{2},\frac{n-p_k}{2})$ provided ${\bf X}$ is of 
full column rank.\\
Let $F$ and $G$ be univariate cumulative distribution functions and $T$ be a test statistic which has distribution function $G$
not necessarily equal to $F$. Let
$p(t|F)=1-F(t)$. By $p$-value of a test statistic $T$ given distribution  $F$ (null distribution) we will mean $p(T|F)$.
We will consider
$p$-values of statistic $R_{jk}^{n}$ given Beta distribution with shape parameters  $\frac{p_k-p_j}{2}$ and $\frac{n-p_k}{2}$. In order to make
notation  simpler $p(R_{jk}^{n}|Beta(\frac{p_k-p_j}{2},\frac{n-p_k}{2}))$ will be denoted as 
$p(R_{jk}^{n}|p_k,p_j)$.
We define the following model selection criteria based on $p$-values of statistic $R_{jk}^{n}$ when one of the indices is held fixed and the other ranges over all potential  models.\\ 
{\bf Minimal $p$-value Criterion (mPVC)}
\[M_m^n={\rm argmin}_{j\in {\cal M}}e^{p_{j}a_{n}} p(R_{0j}^n|p_j,0),\]
 where $p(R_{00}^n|0,0)=e^{a_n}/\sqrt{n}$ and $(a_n)$ is a sequence of nonnegative numbers.
When a minimizer is not unique, the set with the smallest number of elements is chosen. In the case of ties, arbitrary minimizer is selected.
 Observe that when $a_n\equiv 0$ then  from among the pairs $\{(H_0,H_j)\}$ we choose a pair for which we are most inclined to
reject $H_0$ and we select the model corresponding to the most convincing alternative hypothesis.
For positive $a_n$ the scaling factor $e^{p_ja_n}$ is interpreted as additional penalization for the complexity of a model.  
\\
Moreover,  Maximal $p$-value Criterion is defined as\\
{\bf Maximal $p$-value Criterion (MPVC)}
\[M_M^n={\rm argmax}_{j\in {\cal M}}e^{-p_{j}a_{n}} p(R_{jf}^n|M_n,p_j),\]
where $p(R_{ff}^n|M_n,M_n)=1$ and $a_n\to\infty$.
 Thus from  among the pairs $\{(H_j,H_1)\}$ we choose a pair for which we are most reluctant to reject $H_0$ in favour of the full model hypothesis.
We stress that the additional assumption $a_n\to\infty$ needed for consistency of MPVC is not required to prove consistency of mPVC. This point  is discussed  further in  Section 3.
Note that in the definition of both criteria the existence of
encompassing model, either from below or from above, is vital for the construction.  The idea of encompassing has been used in
Bayesian model selection (see e.g. \citet{Casellaetal2009}).\\
Observe that for a fixed number of variables $p_j$  $p$-value $p(R_{0j}^n|p_j,0)$ is a strictly decreasing  function of $R_{0j}^n$. 
Thus the set $M_m^n$ is actually chosen from among subsets for which $R_{0j}^n$ is maximal for the stratum $p_j=1,\ldots,M_n$. The same observation
also holds for MPVC as well as for BIC and AIC. Observe also that if  these criteria choose subsets  of the same cardinality, these subsets necessarily coincide.



\section{Results}
\subsection{Random-design regression} 
\noindent
The main result of this section is consistency of the introduced selectors.
Depending on the context we will use some of the following additional conditions on  the horizons $M_n$, norming constants $a_n$ and matrix $\X$.
\begin{enumerate}
\item[(A1.1')] $M_n/(a_n+\log(n))\to 0$ as $n\to\infty$.
\item[(A1.1'')] $M_n/a_n\to 0$ as $n\to\infty$.
\item[(A1.2)] ${\rm lim}_{n\to\infty} M_n\geq \max_{i\in t} i=:i_{\max}.$
\item[(A1.3)] The minimal eigenvalue $\kappa_n$ of $\E[{\bx_1}^{(n)}{\bx_1}^{(n)'}]$ is bounded away from zero, i.e. $\kappa_n>\kappa>0$ for some $\kappa>0$ and 
$n\in \N$.
\item[(A1.4)] For some $\eta>0$, $n^{-1}{ M}_{n}^{1+\eta}\to 0$ and
\begin{equation}
\label{MomCond}
\sup_{n}\sup_{||{\bf d}||=1}\E|{\bf{d}}'{{\bf z}^{(n)}}|^{4\lceil 2/\eta\rceil}<\infty,
\end{equation}		
where  ${\bf z}^{(n)}=\E[{\bx}_1^{(n)}{\bx}_1^{(n)'}]^{-1/2}\bx_{1}^{(n)}$ is the standardised vector $\bx_{1}^{(n)}$ i.e. 
$E({\bf z}^{(n)}{\bf z}^{(n)'})={\bf I}$ and
$\lceil 2/\eta\rceil$ is the smallest integer greater than or equal to $2/\eta$.
\item[(A1.5)] $a_n/n\to 0$ as $n\to\infty$.

 \end{enumerate}
Assumptions (A1.1') and (A1.1'') are two variants of the condition on a  rate of divergence of $M_n$.  As $M_n$
is nondecreasing, the limit in (A1.2) exists and is either finite or equal to infinity. 
  Condition (A1.2) 
 is a natural condition stating that ultimately the list will contain the true model.
The assumptions (A1.3) and the second part of (A1.4), 
used in \citet{ZhengLoh1997}, imply in particular that with probability tending to one $(\X '\X)^{-1}$ exists and therefore $\hat{\betab}$ is unique. 
Similar conditions are used by \citet{Mammen1993} to study the asymptotic behaviour of bootstrap estimators of  contrasts in 
linear models of increasing dimension.\\
We will consider in detail the case  when $M_{m}^{n}$ and $M_{M}^{n}$ are optimised over all subsets of $f$ i.e. ${\cal M}=2^f$ 
and comment on  the situation  when the nested list of models is considered: 
${\cal M}_{nested}=\{\{1,2,\ldots,i\}\}_{i=1,\ldots,M_n}$.
The first result concerns consistency of the minimal $p$-value criterion.
\begin{Theorem}
\label{Theorem1}
 Let ${\cal M}=2^f$. Then under conditions (A0), (A1.1'), (A1.2), (A1.3), (A1.4), (A1.5)\\
$P(M_{m}^{n}=t)\to 1$, as $n\to\infty$.
\end{Theorem}
As it follows from the  proof an Lemma \ref{LemmaSupmodel} condition (A1.1') may be weakened in  Theorem 1 to $(a_n+\log n-M_n)/\sqrt{M_n}\to\infty$. 
We state now analogous result for MPVC criterion.
\begin{Theorem}
\label{Theorem2}
Let ${\cal M}=2^f$. Then under  conditions of Theorem \ref{Theorem1} with (A1.1') replaced by (A1.1'')\\ $P(M_{M}^{n}=t)\to 1$, as $n\to\infty$.
\end{Theorem}
In order to compare assumptions of the above results note
that when $M_n$ grows  more slowly than $\log(n)$ we can take $a_n=0$ in the case of criterion $M_{m}^{n}$.
However, in the case of $M_{M}^{n}$ the assumption (A1.1'') is obviously not  satisfied for $a_n=0$.\\  
It follows from the proof that the condition (A1.1") may be weakened  in  Theorem \ref{Theorem1} to $(a_n-M_n)/\sqrt{M_n}\to\infty$. \\
Proofs of Theorems \ref{Theorem1} and \ref{Theorem2}  are given in Section 4.
\\Consider now the case when the criteria are optimised over nested list of models
${\cal M}_{nested}=\{\{1,2,\ldots,i\}\}_{i=1,\ldots,M_n}$ 
and define $i_{\max}=\max_{i\in t}i$ as  the largest index of nonzero coefficient in the true model. 
In this case our goal is not to  identify consistently the minimal true model $t$ but rather
$i_{\max}$, which is equivalent to consistent selection of a 
set $t_{\max}=\{1,\ldots,i_{\max}\}$. 
It turns out that  this property holds under weaker conditions than in Theorem \ref{Theorem1} and \ref{Theorem2}.
Namely, the conditions (A1.3) and (A1.4) can be omitted. 
In this case the condition (A0) will be slightly modified. 
Let $\bx_{lt_{\max}}=(x_{l,1}^{(n)},\ldots, x_{l,i_{\max}}^{(n)})'$ be a vector of variables which pertain to the 
model $\{1,\ldots,i_{\max}\}$.
Instead of (A0) we  assume (B0):  $\E(\bx_{lt_{\max}}\bx_{lt_{\max}}')$ is positive definite matrix.  Then under conditions (B0), (A1.1'), (A1.2)and (A1.5) $P(M_{m}^{n}=t_{\max})\to 1$ and analogous result holds for $M_{m}^{n}$ provided (A1.1') is replaced by (A1.1''). This is proved along the lines of the proofs of Theorems \ref{Theorem1} and \ref{Theorem2}.
\\\\In order to lessen computational burden of all subset search we propose two-step model selection with the first step consisting in 
initial ordering of variables according to $p$-values of coefficient of
partial determination (\ref{R}). This method is analogous to the procedure proposed in \citet{ZhengLoh1997} in which variables are ordered according to absolute values of $t$-statistics corresponding to respective attributes.
Then  in the second step an arbitrary criterion Crit is optimised over nested family of models. Specifically, the  greedy procedure consists
of  the following steps. Let
\begin{equation}
\label{greedy}
PV_i=p(R_{(f-\{i\})f}^{n}|M_n,M_n-1),\quad i=1,\ldots,M_n 
\end{equation}
be the $p$-value of statistic $R_{(f-\{i\})f}^{n}$  for testing $H_0:$ model $f-\{i\}$ against $H_1:$ model $f$. 
Then
\begin{enumerate}
\item[(Step 1)] Order the $p$-values in nondecreasing order  $PV_{i_1}\leq PV_{i_2}\leq\ldots\leq PV_{i_{M_n}}$.
\item[(Step 2)]  Consider the nested family $\{\{i_1,i_2,\ldots,i_k\}\}_{k=1,\ldots,M_n}$ and optimise criterion Crit  over this family.
\end{enumerate}
It can be shown that under (A1.2)-(A1.4) 
\begin{equation*}
\lim_{n\to\infty}P(\max_{i\in t}PV_i<\min_{i\not\in t}PV_i)=1.
\end{equation*}
The proof of the above assertion is a simple consequence of Theorem 2 in \citet{ZhengLoh1997}.
This, together with Theorems \ref{Theorem1} and \ref{Theorem1} for the case of the nested list of models, when minimal or maximal $p$-value 
criterion is considered as Crit,  leads to the following corollary.
\begin{Corollary}
Under conditions of Theorems \ref{Theorem1} and \ref{Theorem2} respectively the greedy versions of mPVC and MPVC procedures are consistent.
\end{Corollary}
Observe that  since parameters of  beta distribution used to calculate $p$-values in (\ref{greedy}) do not change with $i$,
the ordering in the first step is equivalent to ordering wrt values of $R_{(f-\{i\})f}^{n}$, or  to the ordering wrt to
absolute values of $t$-statistics when the full model is fitted.\\ 

\subsection{Deterministic-design regression}
In this section we will briefly discuss the case when the design matrix $\X$ is nonrandom. We allow that the values of attributes $\bx_{l,1}^{(n)},\ldots,\bx_{l,M_n}^{(n)}$ of $l^{\textrm{th}}$ observation may depend on $n$.
Recall that $\bx_{lt}=\bx_{lt}^{(n)}$ is a vector of variables which pertain to the minimal true model $t$.
In the case of all subset search we  replace condition (A0) by the following assumption
\begin{enumerate}
\item[(C0)] $n^{-1}\sum_{i=1}^{n}\bx_{lt}\bx_{lt}'\to \Wa$, as $n\to\infty$, where $\Wa$ is a positive definite matrix. 
\end{enumerate}
In the case of random covariates the above convergence in probability  follows from The Law of Large Numbers.
We  also replace conditions (A1.3) and (A1.4) by the following assumption
\begin{enumerate}
\item[(C1)] The minimum eigenvalue $\tilde{\kappa}_n$ of $n^{-1}\X '\X$ is bounded away from zero, i.e. $\tilde{\kappa}_n>\tilde{\kappa}>0$ for some $\tilde{\kappa}>0$ and $n\in \N$.
\end{enumerate}
Recall that $\hat{\betab}=(\hat{\beta}_1,\ldots,\hat{\beta}_{M_n})'$ is the least squares estimator based on the full model $f$. Let $T_i=\hat{\sigma}^{-1}[(\X '\X)^{-1}_{i,i}]^{-1/2}$ be the corresponding t-statistic. It can be easily shown that $\hat{\sigma}T_i=\beta_i[(\X '\X)^{-1}_{i,i}]^{-1/2}+o_{P}(1)$, for $i\in t$. Thus by assumption (C1) $P(\hat{\sigma}T_i>Cn^{-1/2})\to 1$ as $n\to\infty$, for some $C>0$. This implies  the conclusion of Lemma \ref{LemmaSubmodel} in Section 4, namely that for $i\in t$ with probability tending to one $RSS(f-\{i\})/RSS(f)$ is bounded away  from $0$.
As (A1.3) and (A1.4) are used in the random-design case only to prove Lemma \ref{LemmaSubmodel} it follows that the analogous results to Theorem \ref{Theorem1} and Theorem \ref{Theorem2} hold for the deterministic-design case.
\begin{Corollary}
\label{Col1}
Under conditions (C0), (A1.1'), (A1.2), (C1), (A1.5) \\
$P(M_{m}^{n}=t)\to 1$, as $n\to\infty$. 
\end{Corollary}
\begin{Corollary}
\label{Col2}
Under conditions of Corollary \ref{Col1} with (A1.1') replaced by (A1.1'')\\
$P(M_{M}^{n}=t)\to 1$, as $n\to\infty$. 
\end{Corollary}
Consider the case of nested family search. Recall that $\bx_{lt_{\max}}$ is a vector of variables which pertain to the model $\{1,\ldots,i_{\max}\}$.
If condition (B0) if replaced by the following assumption
\begin{enumerate}
\item[(D0)] $n^{-1}\sum_{i=1}^{n}\bx_{lt_{\max}}\bx_{lt_{\max}}'\to \Wb$, as $n\to\infty$, where $\Wb$ is a positive definite matrix. 
\end{enumerate}
then results discussed at the end of Section 3.1 hold for deterministic design.
\section{Proofs}
We first state auxiliary lemmas which will be used in the proof of Theorem \ref{Theorem1}. The first one proved in \citet{PokarowskiMielniczuk2010} gives
an approximation of tail probability function of beta distribution. 
Let $B_{a,b}$ be a random variable having  beta distribution with shape parameters $a$ and $b$ and $B(x,y)$ denote beta function. 
Define an 
auxiliary function
\begin{equation*}
L(a,b,x)=\frac{(a-1)(1-x)}{1-a+(a+b)x},
\end{equation*}
for $a,b,x\in\R$ such that $x\neq (a-1)/(a+b)$.
\begin{Lemma}
\label{LemmaBeta}
Assume $x>\frac{a-1}{a+b}$. Then  for $a\geq 1$ 
\begin{equation}
\label{LemmaBetaeq1}
\frac{(1-x)^{b}x^{a-1}}{B(a,b)b}\leq P[B_{a,b}>x]\leq
\frac{(1-x)^{b}x^{a-1}}{B(a,b)b}(1+L(a,b,x)) 
\end{equation}
and   for $a < 1$ 
\begin{equation}
\label{LemmaBetaeq2}
\frac{(1-x)^{b}x^{a-1}}{B(a,b)b}(1+L(a,b,x))\leq
P[B_{a,b}>x]\leq
\frac{(1-x)^{b}x^{a-1}}{B(a,b)b} .
\end{equation}
\end{Lemma}
The following Lemma states simple but useful inequalities for gamma function.
\begin{Lemma}
\label{LemmaGamma}
Let $a=p/2$ and $b=(n-p)/2$, for some $p,n\in\N$. Then
\begin{equation*}
\Gamma(b)b^{a}\leq\Gamma(a+b)\leq\frac{2}{\sqrt{\pi}}\Gamma(b)(a+b)^a.
\end{equation*}
\end{Lemma}
The above Lemma implies an inequality for  beta function $B(a,b)=\Gamma(a)\Gamma(b)/\Gamma(a,b)$
\begin{equation}
\label{LemmaGammaEq}
\frac{b^{a-1}}{\Gamma(a)}\leq\frac{1}{bB(a,b)}\leq\frac{2}{\sqrt{\pi}}\frac{(a+b)^a}{b\Gamma(a)},
\end{equation}
for $a=p/2$, $b=(n-p)/2$ and $p,n\in\N$.\\
\begin{Remark}
Lemma \ref{LemmaGamma} easily implies inequality $\Gamma(p/2)\leq (\lceil p/2\rceil-1)!\leq p^{p/2}$ for $p>1$, which will be frequently used throughout.
\end{Remark}
The following Lemma states that for a proper submodel of the true model $t$ variance estimator is asymptotically biased. $j\subset k$ denotes a proper inclusion of $j$ in $k$.
\begin{Lemma}
\label{Auxilliary}
(i)  For $j\supseteq t$, $j\in{\cal M}$ 
$\frac{RSS(j)}{n}\cp\sigma^{2}$ as $n\to\infty$.
Moreover, for  $j\subset t$, $j\in{\cal M}$ if  (A0) is satisfied then  $\frac{RSS(j)}{n}\cp\sigma^{2}+\lambda_{j}$ as $n\to\infty$, where 
$\lambda_{j}>0$ .\\
(ii) Let $j\subset t_{\max}$, $j\in {\cal M}_{nested}$ and assume (B0). Then $\frac{RSS(j)}{n}\cp\sigma^{2}+\lambda_{j}$ as $n\to\infty$, where $\lambda_{j}>0$ .\\
\end{Lemma} 
\begin{Lemma}
\label{LemmaSupmodel}
Let $R_n$ be a sequence of real numbers such that 
$(R_n-M_n)/\sqrt{M_n}\to\infty$ as $n\to\infty$. 
Assume also that $M_n/n\to 0$ and matrix $\X '\X$ is invertible with probability tending to 1. Then 
\begin{equation*}
P\left\{n\log\left[\frac{RSS(t)}{RSS(f)}\right]>R_n\right\}\to 0 
\end{equation*}
as $n\to\infty$.
\end{Lemma} 
\begin{Remark}
Observe that as $(R_n-M_n)/\sqrt{M_n}=\sqrt{M_n}(R_n/M_n-1)$, the imposed condition on $R_n$ is implied by $R_n/M_n\to\infty$.
Thus in particular 
Lemma \ref{LemmaSupmodel} implies that 
\begin{equation*}
\frac{RSS(t)}{RSS(f)}=O_{P}\left[\exp\left(\frac{R_n}{n}\right)\right],
\end{equation*}
for any $R_n$ such that $R_n/M_n\to\infty$.
Observe moreover that Lemma \ref{LemmaSupmodel} holds true also in the case  $M_n=M$  when the condition on $R_n$ reduces to $R_n\to\infty$
only and thus $RSS(t)/RSS(f)=O_P(exp(n^{-1}))$. This can be seen directly 
 from Lemma \ref{Auxilliary} and the fact that $RSS(t)-RSS(f)\sim \chi^2_{M-p_t}$ as  it follows from them that  $R_{tf}^n={\cal O}_P(n^{-1})$ and thus
$n\log(RSS(t)/RSS(f))={\cal O}_P(1)$. 
\end{Remark}
\begin{Lemma}
\label{LemmaSubmodel}
Assume conditions (A1.3) and (A1.4). 
Then there exists $a>0$ such that
\begin{equation*}
P\left\{\min_{i\in t} \log\left[\frac{RSS(f-\{i\})}{RSS(f)}\right]>a\right\}\to 1 
\end{equation*}
as $n\to\infty$.
\end{Lemma} 
Thus Lemma \ref{LemmaSubmodel} implies that with probability tending to $1$ 
$RSS(f-\{i\})/{RSS(f)}$ for $i\in t$ is bounded away from 0.

\subsection{Proof of Theorem \ref{Theorem1}}
We will consider separately two  cases: the first when the true model $t$ contains nontrivial regressors ($p_t\geq 1$) and
the second, when it equals the null model.\\
{\bf{Case 1} ($p_t\geq 1$).} We will treat the case $p_t\geq 2$ in detail, the case $p_t=1$ is similar but simpler and relies
on (\ref{LemmaBetaeq2}) instead of (\ref{LemmaBetaeq1}) to treat $\mPVC{t}$.\\
(i) Let  $j$ be such that $j\supset t$ i.e. $t$ is a proper subset of $j$. 
We will prove that $P[e^{p_ta_n}\mPVC{t}>\inf_{j\supset t}e^{p_ja_n}\mPVC{j}]\to 0$ as $n\to\infty$.  
Using (\ref{LemmaGammaEq})  with $a=p_t/2$ and $b=(n-p_t)/2$ we obtain the following inequalities  for sufficiently large $n$  
\begin{equation}
\label{D2}
\frac{1}{B(\frac{p_t}{2},\frac{n-p_t}{2})\left(\frac{n-p_t}{2}\right)}\leq
\frac{2\left(\frac{n}{2}\right)^{\frac{p_t}{2}}}{\sqrt{\pi}\left(\frac{n-p_t}{2}\right)\Gamma\left(\frac{p_t}{2}\right)}\leq
\frac{2\left(\frac{n}{2}\right)^{\frac{p_t}{2}}}{\sqrt{\pi}\left(\frac{n}{4}\right)\Gamma\left(\frac{p_t}{2}\right)}=
\frac{4\left(\frac{n}{2}\right)^{\frac{p_t}{2}-1}}{\sqrt{\pi}\Gamma\left(\frac{p_t}{2}\right)}.
\end{equation}
Moreover for $j\supset t$ and sufficiently large $n$
\begin{equation}
\label{D3}
\frac{1}{B(\frac{p_j}{2},\frac{n-p_j}{2})\left(\frac{n-p_j}{2}\right)}\geq
\frac{\left(\frac{n-p_j}{2}\right)^{\frac{p_j}{2}-1}}{\Gamma\left(\frac{p_j}{2}\right)}\geq
\frac{\left(\frac{n-M_n}{2}\right)^{\frac{p_j}{2}-1}}{M_n^{\frac{p_j}{2}}}\geq
\frac{\left(\frac{n-M_n}{2}\right)^{\frac{p_t+1}{2}-1}}{M_n^{\frac{p_t+1}{2}}}\geq
\frac{\left(\frac{n}{2}\right)^{\frac{p_t+1}{2}-1}\left(\frac{1}{2}\right)^{\frac{p_t+1}{2}-1}}{M_n^{\frac{p_t+1}{2}}}.
\end{equation}
Note that 
\begin{eqnarray*}
&&
P\left(\inf_{j\supset t}R_{0j}^{n}\geq \sup_{j\supset t}\frac{\frac{p_j}{2}-1}{\frac{n}{2}}\right)\leq
P\left(R_{0t}^{n}\geq (M_n-2)/n\right)\to 1,
\end{eqnarray*}
which follows from Lemma \ref{Auxilliary} and the fact that $M_n/n\to 0$. 
Thus the assumption of Lemma \ref{LemmaBeta} is satisfied for $x=R_{0j}^{n}$, $a=\frac{p_j}{2}$, $b=\frac{n-p_j}{2}$ and all $j\supset t$.
Using (\ref{LemmaBetaeq1}) we have
\begin{eqnarray}
\label{D4}
&&
P[e^{p_ta_n}\mPVC{t}>\inf_{j\supset t}e^{p_ta_n}\mPVC{j}]\leq
\cr
&&
P\left\{\frac{(1-R_{0t}^{n})^{\frac{n-p_t}{2}}(R_{0t}^{n})^{\frac{p_t}{2}-1}[1+L\left(\frac{p_t}{2},\frac{n-p_t}{2},R_{0t}^{n}\right)]e^{p_ta_n}}{B(\frac{p_t}{2},\frac{n-p_t}{2})\left(\frac{n-p_t}{2}\right)}>
\inf_{j\supset t}
\frac{(1-R_{0j}^{n})^{\frac{n-p_j}{2}}(R_{0j}^{n})^{\frac{p_j}{2}-1}e^{(p_t+1)a_n}}{B(\frac{p_j}{2},\frac{n-p_j}{2})\left(\frac{n-p_j}{2}\right)}
\right\}\leq
\cr
&&
P\left\{\frac{(1-R_{0t}^{n})^{\frac{n-p_t}{2}}[1+L\left(\frac{p_t}{2},\frac{n-p_t}{2},R_{0t}^{n}\right)]e^{p_ta_n}}{B(\frac{p_t}{2},\frac{n-p_t}{2})\left(\frac{n-p_t}{2}\right)}>
\inf_{j\supset t}
\frac{(1-R_{0f}^{n})^{\frac{n-p_t}{2}}(R_{0t}^{n})^{\frac{M_n}{2}-1}e^{(p_t+1)a_n}}{B(\frac{p_j}{2},\frac{n-p_j}{2})\left(\frac{n-p_j}{2}\right)}
\right\}.
\end{eqnarray}
Taking logarithms and using inequalities (\ref{D2}), (\ref{D3}) we obtain
\begin{eqnarray*}
&&
P[\log\mPVC{t}+p_ta_n>\inf_{j\supset t}\log\mPVC{j}+(p_t+1)a_n]\leq
P\left\{\left[\frac{n-p_t}{2}\right]\log\left[\frac{RSS(t)}{RSS(f)}\right]>\tilde{W}_n\right\},
\end{eqnarray*}
where
\begin{eqnarray*}
&&
\tilde{W}_n=a_n+\frac{1}{2}\log\left(\frac{n}{2}\right)-\log[1+L\left(\frac{p_t}{2},\frac{n-p_t}{2},R_{0t}^{n}\right)]+
\left(\frac{M_n}{2}-1\right)\log(R_{0t}^n)+
\cr
&&
\left(\frac{p_t+1}{2}-1\right)\log\left(\frac{1}{2}\right)-\left(\frac{p_t+1}{2}\right)\log(M_n)-\log\left(\frac{4}{\sqrt{\pi}}\right)+\log\Gamma\left(\frac{p_t}{2}\right).
\end{eqnarray*}
Assumption $M_n/(a_n+\log(n))\to 0$, Lemma \ref{Auxilliary} and the fact that $R_{0,t}\cp\sigma^2>0$
 imply that there exists a sequence $W_n$ of real numbers such that $P(\tilde{W}_n>W_n)\to 1$ and $W_n/M_n\to\infty$. Now the required convergence follows from
\begin{equation*}
P\left\{\left[\frac{n-p_t}{2}\right]\log\left[\frac{RSS(t)}{RSS(f)}\right]>W_n\right\}\to 0
\end{equation*}
which in its turn is implied by Lemma \ref{LemmaSupmodel}.\\ 
(ii) Consider now the case $j\nsupseteq t$ and let $i=i(j)\in\N$ be such that $i\in t\cap j^c$.
 We will prove that  $P[e^{p_ta_n}\mPVC{t}>\inf_{j\nsupseteq t}e^{p_ja_n}\mPVC{j}]\to 0$ as $n\to\infty$. 
Define $M(n,i)=\max\{R_{0(f-\{i\})}^{n},\frac{2M_n}{(n-M_n)}\}$, for $i\in t$. 
Assume first that $p_j\geq 2$.
Using (\ref{LemmaBetaeq1}) and (\ref{LemmaGammaEq}) we have
\begin{eqnarray}
\label{D8}
&&
e^{p_ja_n}p(R_{0j}^{n}|p_j,0)\geq e^{2a_n}p(M(n,i)|p_j,0)\geq
\frac{e^{2a_n}[1-M(n,i)]^{\frac{n-p_j}{2}}M(n,i)^{\frac{p_j}{2}-1}}{B\left(\frac{p_j}{2},\frac{n-p_j}{2}\right)\left(\frac{n-p_j}{2}\right)}\geq
\cr
&&
\frac{e^{2a_n}[1-M(n,i)]^{\frac{n}{2}}\left(\frac{2M_n}{n-M_n}\right)^{\frac{p_j}{2}-1}}{B\left(\frac{p_j}{2},\frac{n-p_j}{2}\right)\left(\frac{n-p_j}{2}\right)}\geq
\frac{e^{2a_n}[1-M(n,i)]^{\frac{n}{2}}\left(\frac{2M_n}{n-M_n}\right)^{\frac{p_j}{2}-1}\left(\frac{n-M_n}{2}\right)^{\frac{p_j}{2}-1}}{\Gamma\left(\frac{p_j}{2}\right)}\geq
\cr
&&
\frac{e^{2a_n}[1-M(n,i)]^{\frac{n}{2}}M_n^{\frac{p_j}{2}-1}}{M_n^{\frac{p_j}{2}}}=e^{2a_n}[1-M(n,i)]^{\frac{n}{2}}M_n^{-1}.
\end{eqnarray}
From (\ref{LemmaBetaeq1}) and (\ref{D2})
\begin{eqnarray}
\label{D9}
e^{p_ta_n}p(R_{0t}^{n}|p_t,0)\leq\frac{e^{p_ta_n}(1-R_{0t}^{n})^{\frac{n-p_t}{2}}4\left(\frac{n}{2}\right)^{\frac{p_t}{2}-1}\left[1+L\left(\frac{p_t}{2},\frac{n-p_t}{2},R_{0t}^{n}\right)\right]}{\sqrt{\pi}\Gamma\left(\frac{p_t}{2}\right)}
\end{eqnarray}
Using (\ref{D8}) and (\ref{D9}) we have for $p_t\geq 2$ and $p_j\geq 2$
\begin{eqnarray*}
&&
P[e^{p_ta_n}\log\mPVC{t}>\inf_{j\nsupseteq t}e^{p_ja_n}\log\mPVC{j}]\leq
P\left\{\inf_{i\in t}\frac{n}{2}\log\left[\frac{(1-M(n,i))RSS(0)}{RSS(t)}\right]<\tilde{S}_n\right\},
\end{eqnarray*}
where
\begin{eqnarray*}
&&
\tilde{S}_n=a_n(p_t-2)+\left(\frac{p_t}{2}-1\right)\log\left(\frac{n}{2}\right)-\frac{p_t}{2}\log\left(\frac{RSS(t)}{RSS(0)}\right)+
\log\left(\frac{4}{\sqrt{\pi}}\right)+
\cr
&&
\log\left[1+L\left(\frac{p_t}{2},\frac{n-p_t}{2},R_{0t}^{n}\right)\right]+
\log\Gamma^{-1}\left(\frac{p_t}{2}\right)+\log(M_n).
\end{eqnarray*}
In view of definition of $M(n,i)$ the last probability can be bounded from above by 
\begin{eqnarray*}
&&
P\left\{\inf_{i\in t}\frac{n}{2}\log\left[\frac{RSS(f-\{i\})}{RSS(t)}\right]<\tilde{S}_n\right\}+
P\left\{\frac{n}{2}\log\left[\frac{(1-\frac{2M_n}{n-M_n})RSS(0)}{RSS(t)}\right]<\tilde{S}_n\right\}.
\end{eqnarray*}
The second probability above converges to zero in view of Lemma \ref{Auxilliary}. Consider the first probability. 
Since the number of elements of $t$ is finite 
it suffices show that $P\left\{\frac{n}{2}\log\left[\frac{RSS(f-\{i\})}{RSS(t)}\right]<\tilde{S}_n\right\}\to 0$ for any $i\in t$. Namely, it 
is bounded  from above by
\begin{eqnarray}
\label{D11}
&&
P\left\{\frac{n}{2}\log\left[\frac{RSS(f-\{i\})}{RSS(f)}\right]+
\frac{n}{2}\log\left[\frac{RSS(f)}{RSS(t)}\right]<\tilde{S}_n\right\}\leq
\cr
&&
P\left\{\frac{n}{2}\log\left[\frac{RSS(f-\{i\})}{RSS(f)}\right]<2\tilde{S}_n\right\}+
P\left\{\frac{n}{2}\log\left[\frac{RSS(f)}{RSS(t)}\right]<-\tilde{S}_n\right\}\leq
\cr
&&
P\left\{n\log\left[\frac{RSS(f-\{i\})}{RSS(f)}\right]<\tilde{S}_n\right\}+
P\left\{\frac{n}{2}\log\left[\frac{RSS(t)}{RSS(f)}\right]\geq\tilde{S}_n\right\}.
\end{eqnarray}
From assumptions (A1.5) and (A1.1') $\tilde{S}_n/n\cp 0$ and $\tilde{S}_n/M_n\cp\infty$, respectively. Thus 
the convergence to zero of the above two probabilities in (\ref{D11}) follows from Lemma
\ref{LemmaSubmodel} and \ref{LemmaSupmodel}, respectively.
The case $p_j=1$ is treated  analogously.\\
Consider now the case $p_j=0$. From (\ref{D9}) we have
\begin{eqnarray}
\label{D12}
&&
P[\log p(R_{0t}^{n}|p_t,0)+p_ta_n>\log p(R_{00}^{n}|0,0)]=
P[\log p(R_{0t}^{n}|p_t,0)>a_n-\frac{1}{2}\log(n)-p_ta_n]\leq
\cr
&&
P\left\{\left(\frac{n-p_t}{2}\right)\log\left[\frac{RSS(0)}{RSS(t)}\right]<G_n\right\},
\end{eqnarray}
where
\begin{eqnarray*}
&&
G_n=(p_t-1)a_n+\frac{1}{2}\log(n)+\left(\frac{p_t}{2}-1\right)\log\left(\frac{n}{2}\right)+\log\left(\frac{4}{\sqrt{\pi}}\right)
+\log\Gamma^{-1}\left(\frac{p_t}{2}\right)+
\cr
&&
\log\left[1+L\left(\frac{p_t}{2},\frac{n-p_t}{2},R_{0t}^{n}\right)\right]
.
\end{eqnarray*}
The convergence to zero of the probability in (\ref{D12}) follows from Lemma \ref{Auxilliary} amd assumption (A1.5).\\
{\bf Case 2 $(p_t=0)$} i.e. the true model is null model. We treat in detail the case $p_j\geq 2$. 
Define $\bar{M}(n)=\max\{R_{0f}^{n},\frac{2M_n}{n-M_n}\}$. Note that 
the assumption of Lemma \ref{LemmaBeta} is satisfied for $x=\bar{M}(n)$, $a=\frac{p_j}{2}$, and $b=\frac{n-p_j}{2}$.
Using (\ref{LemmaBetaeq1}) and (\ref{LemmaGammaEq}) we have
\begin{eqnarray}
\label{D5}
&&
e^{p_ja_n}p(R_{0j}^{n}|p_j,0)\geq e^{2a_n}p(\bar{M}(n)|p_j,0)\geq
\frac{e^{2a_n}[1-\bar{M}(n)]^{\frac{n-p_j}{2}}\bar{M}(n)^{\frac{p_j}{2}-1}}{B\left(\frac{p_j}{2},\frac{n-p_j}{2}\right)\left(\frac{n-p_j}{2}\right)}\geq
\cr
&&
\frac{e^{2a_n}[1-\bar{M}(n)]^{\frac{n-p_t}{2}}\left(\frac{2M_n}{n-M_n}\right)^{\frac{p_j}{2}-1}}{B\left(\frac{p_j}{2},\frac{n-p_j}{2}\right)\left(\frac{n-p_j}{2}\right)}\geq
\frac{e^{2a_n}[1-\bar{M}(n)]^{\frac{n-p_t}{2}}\left(\frac{2M_n}{n-M_n}\right)^{\frac{p_j}{2}-1}\left(\frac{n-M_n}{2}\right)^{\frac{p_j}{2}-1}}{\Gamma\left(\frac{p_j}{2}\right)}\geq
\cr
&&
\frac{e^{2a_n}[1-\bar{M}(n)]^{\frac{n-p_t}{2}}M_n^{\frac{p_j}{2}-1}}{M_n^{\frac{p_j}{2}}}=e^{2a_n}[1-\bar{M}(n)]^{\frac{n-p_t}{2}}M_n^{-1}.
\end{eqnarray}
Using (\ref{D5}) we obtain the following inequality 
\begin{eqnarray}
\label{D6}
&&
P[\log p(R_{00}^{n}|0,0)>\inf_{j:p_j\geq 2}\log\mPVC{j}+2a_n]\leq
P[a_n-\frac{1}{2}\log(n)>\inf_{j:p_j\geq 2}\log\mPVC{j}+2a_n]\leq
\cr
&&
P\left\{-\left(\frac{n-p_t}{2}\right)\log[1-\bar{M}(n)]>a_n+\frac{1}{2}\log(n)-\log(M_n)\right\}\leq
\cr
&&
P\left\{\left(\frac{n-p_t}{2}\right)\log\left[\frac{RSS(0)}{RSS(f)}\right]>a_n+\frac{1}{2}\log(n)-\log(M_n)\right\}+
\cr
&&
I\left[-\left(\frac{n-p_t}{2}\right)\log\left(1-\frac{2M_n}{n-M_n}\right)>a_n+\frac{1}{2}\log(n)-\log(M_n)\right].
\end{eqnarray}
From Lemma \ref{LemmaSupmodel} and the assumption $M_n/(a_n+\log(n))\to 0$ the first probability in (\ref{D6}) converges to zero.
The same assumption implies that the second term is ultimately 0. 
This completes the proof.\\

\subsection{Proof of Theorem \ref{Theorem2}}
The proof is similar to that of Theorem 1 and splits into  two cases: $M_n-p_t\geq 1$ 
(corresponding to the case $p_t\geq 1$ in the previous proof) and $M_n=p_t$ (corresponding to the former 
case $p_t=0$). We  give the sketch of the proof only.
\\ 
{\bf Case 1 $(M_n- p_t\geq 1)$}. We discuss the situation when $M_n- p_t\geq 2$, the remaining case relies on (\ref{LemmaBetaeq2}) instead of (\ref{LemmaBetaeq1}).
Define $\tilde{M}(n,t)=\max\{R_{tf}^{n},\frac{2M_n}{n-M_n}\}$. Note that 
the assumption of Lemma \ref{LemmaBeta} is satisfied for $x=\tilde{M}(n,t)$, $a=\frac{M_n-p_t}{2}$, and $b=\frac{n-M_n}{2}$.
In this case condition $a\geq 1$ is also satisfied.
Analogously to the proof of (\ref{D5}) we obtain 
\begin{equation}
\label{DD1}
p(R_{tf}^{n}|M_n,p_t)\geq p(\tilde{M}(n,t)|M_n,p_t)\geq [1-\tilde{M}(n,t)]^{\frac{n-M_n}{2}}M_n^{-1}.
\end{equation}
(i) Let  $j$ be such that $j\supset t$ i.e. $t$ is a proper subset of $j$. 
We will prove that\\
$P[\MPVC{t}<\sup_{j\supset t}\MPVC{j}]\to 0$ as $n\to\infty$.  
For $j\supset t$ we have $\MPVC{j}\leq\exp[-(p_t+1)a_n]$. This inequality also applies to $j=f$. Thus using (\ref{DD1}) we 
obtain the following inequalities
\begin{eqnarray*}
&&
P[\MPVC{t}<\sup_{j\supset t}\MPVC{j}]\leq
\cr
&&
P\left\{\left(\frac{n-M_n}{2}\right)\log[1-\tilde{M}(n,t)]-\log(M_n)-p_ta_n<-(p_t+1)a_n\right\}\leq
\cr
&&
P\left\{\left(\frac{n-M_n}{2}\right)\log\left[\frac{RSS(t)}{RSS(f)}\right]>a_n-\log(M_n)\right\}+
\cr
&&
I\left\{-\left(\frac{n-M_n}{2}\right)\log\left[1-\frac{2M_n}{n-M_n}\right]>a_n-\log(M_n)\right\}.
\end{eqnarray*}
The above bound converges  to zero in view of the assumption $M_n/a_n\to 0$ and Lemma \ref{LemmaSupmodel}.\\
(ii) Consider now the case $j\nsupseteq t$ and assume that $p_j\leq M_n-2$ (this corresponds to $p_j\geq 2$ in the previous proof).
Let index $i=i(j)$ be such that $i\in t\cap j^c$. 
It follows from Lemma \ref{LemmaSubmodel} that
the assumption of Lemma \ref{LemmaBeta} is satisfied for $x=R_{(f-\{i\})f}$, $a=\frac{M_n-p_j}{2}$, and $b=\frac{n-M_n}{2}$.
Moreover the same reasoning yields for all $j\nsupseteq t$ 
$L\left(\frac{M_n-p_j}{2},\frac{n-M_n}{2},R_{f-\{i\}f}\right)\leq M_n$ wih probability tending to 1.
Using (\ref{LemmaBetaeq1})  we have the following inequalities
\begin{eqnarray}
\label{DD4}
&&
\MPVC{j}\leq p(R_{(f-\{i\})f}|M_n,p_j)\leq
\cr
&&
\frac{[1-R_{(f-\{i\})f}]^{\frac{n-M_n}{2}}[R_{(f-\{i\})f}]^{\frac{M_n-p_j}{2}-1}}
{B\left(\frac{M_n-p_j}{2},\frac{n-M_n}{2}\right)\left(\frac{n-M_n}{2}\right)}
\left[1+L\left(\frac{M_n-p_j}{2},\frac{n-M_n}{2},R_{(f-\{i\})f}\right)\right]\leq
\cr
&&
[1-R_{(f-\{i\})f}]^{\frac{n-M_n}{2}}\frac{2n^{\frac{M_n}{2}}}{\sqrt{\pi}\Gamma\left(\frac{M_n}{2}\right)}[1+M_n].
\end{eqnarray}

Thus
\begin{eqnarray*}
&&
P[\MPVC{t}<\sup_{j\nsupseteq t}\MPVC{j}]\leq
\cr
&&
P\left\{\sup_{i\in t}\left(\frac{n-M_n}{2}\right)\log\left[\frac{\left(1-\tilde{M}(n,t)\right)RSS(f-\{i\})}{RSS(f)}\right]<K_n\right\},
\end{eqnarray*}
where
\begin{eqnarray*}
&&
K_n=p_ta_n+\log\left(\frac{2}{\sqrt{\pi}}\right)+\frac{M_n}{2}\log(n)-
\log\Gamma\left(\frac{M_n}{2}\right)+\log(1+M_n)+\log(M_n).
\end{eqnarray*}
Similarly to the proof of (\ref{D11})  we obtain that the RHS tends to 0.
\\
The case $p_j>M_n-2$ is simpler and uses (\ref{LemmaBetaeq2}) instead of (\ref{LemmaBetaeq1}).\\
{\bf Case 2 $(M_n=p_t)$}. 
Thus $\MPVC{t}=e^{-M_na_n}$. Assume $p_j\leq M_n-2$ and let $i=i(j)$ be such that $i\in j^c\cap t$.  Then  using 
$L\left(\frac{M_n-p_j}{2},\frac{n-M_n}{2},R_{jf}\right)\leq M_n$ and (\ref{LemmaBetaeq1}) (cf (\ref{DD4})
it is easy to establish
that
\[\MPVC{j}\leq p(R_{(f-\{i\})f}|M_n,p_j)\leq
[1-R_{(f-\{i\})f}]^{\frac{n-M_n}{2}}\frac{2n^{\frac{M_n}{2}}}{\sqrt{\pi}\Gamma\left(\frac{M_n}{2}\right)}[1+M_n].\]
Then it follows that
\begin{eqnarray*}
&&
P[\MPVC{t}<\sup_{j\not = t}\MPVC{j}]\leq
\cr
&&
P\left\{\sup_{i\in t}\left(\frac{n-M_n}{2}\right)\log\left[\frac{RSS(t-\{i\})}{RSS(t)}\right]<\tilde{K}_n\right\},
\end{eqnarray*}
where
\begin{eqnarray*}
&&
\tilde{K}_n=M_na_n+\log({2}/\sqrt{\pi})+\frac{M_n}{2}\log(n)-
\log\Gamma(M_n/2)+\log(1+M_n)+\log(M_n).
\end{eqnarray*}
The convergence to zero of the above probability follows from Lemma \ref{Auxilliary} and the assumption $a_n/n\to 0$.  The case 
$p_j>M_n-2$ is analogous.


\section{Numerical experiments}
In this section we study the finite-sample performance of the model selection procedures. We consider criteria defined in Section 2: minimal p-value criterion $M^{n}_{m}$ with $a_n=0$ which will be called simply in this section mPVC and two scaled p-value criteria with scalings which were empirically chosen, namely minimal p-value criterion with $a_n=\log(n)/2$  and maximal p-value criterion with the same $a_n$ called mPVCcal and MPVCcal, respectively.
As benchmarks  we  considered  performance  of classical criteria based on penalized log-likelihood which have the form
\begin{equation*}
\label{BIC} 
{\rm argmax}_{j\in {\cal M}}\{2\log{\bf f}_{\hat{\betab_j},{\hat\sigma_j}^2}(\Y|\X)-p_jC_n\}
={\rm argmax}_{j\in {\cal M}}\{-n\log[RSS(j)/n] -p_jC_n\}
\end{equation*}
with penalties: $C_n=2$ and $C_n=\log(n)$ which correspond to Akaike (AIC) and Bayesian (BIC) information criteria, respectively.
\subsection{Simulation experiments}
The simulation experiments were carried out with sample sizes $n=75, 100, 200, 300, 500, 1000$ repeated $N=500$ times.   
We consider the following lists of  models
\begin{enumerate}
\item[(M1)] $t=\{10\}$, $\beta_1=0.2$, $M_n=30$,
\item[(M2)] $t=\{1,2,5,6\}$, $\betab=(0.9,-0.8,-0.4,0.2)'$, $M_n=6$,
\item[(M3)] $t=\{2,4,5\}$, $\betab=(1,1,1)'$, $M_n=5$,
\item[(M4)] $t=\{2k+7:k=3,\ldots,12\}$, $\betab=(1,\ldots,1)'$, $M_n=60$.
 \end{enumerate}
In all cases ${\cal M}=2^{\{1,\ldots,M_n\}}$.
Models M1, M3 and M4 were also considered in \citet{ZhengLoh1997}. 
Regressors $\bx^{n}_{l}$ were generated from
$M_n$-variate zero mean normal distribution with $(i,j)$th entry of the covariance matrix $\Sigma_{X}=(\sigma_{ij})_{ij}$ equal 
$\sigma_{ij}=0.5^{|i-j|}$. The distribution of $(\varepsilon_1,\ldots,\varepsilon_n)$ was multivariate standard normal.
We considered greedy variants of the selection methods,
described in Section 3. Table \ref{TABORDERING} presents estimated probabilities of correct ordering, e.g. the probabilities that the coordinates corresponding to nonzero coefficients are placed ahead the spurious ones. It is seen that for $n\geq 500$ for the models considered a correct ordering is recovered practically always. 
We assess the effectiveness of the selection rule in terms of the probability of true model selection $P(\hat{t}=t)$, 
where $\hat{t}$ is a model selected by the considered rule and mean squared error $\E(||\X\betab-\X\hat{\betab}(\hat{t})||^{2})$, 
where $\hat{\betab}(\hat{t})$ is the post-model selection estimator of $\betab$ i.e. ML estimator in the chosen model.
In the experiments estimates of these
measures calculated as the empirical means of respective quantities were considered. 
The influence of the sample size on the effectiveness of selected rules has been investigated. For models M1, M3 and M4 criterion MPVCcal and mPVCcal perform considerably better for all sample sizes considered than mPVC and commonly used BIC and AIC (see Figure 1 and 2).
In contrast, in the case of model M2 criterion mPVC works better than others. 
In general, performance of mPVCcal is similar to that of MPVCcal.
The results also indicate that model M1 with the only one significant variable placed at position 10 is the most difficult for selection among the models considered.
This is due to the fact that  in this case it is difficult to recover the correct ordering (see Table \ref{TABORDERING}), especially for small sample sizes. Secondly the selection criteria seem to work worse when the number of nuisance covariates is large. 
For model M1 we also studied the influence of the value of the true parameter $\beta_1$. Figure 3 indicates that performance of both measures is much worse for small values of the parameter.
The influence of the size of the list $M_n$ on the effectiveness of selection rules has been also investigated. Figure 4 shows that for model M1 performance of the AIC, BIC and mPVC is influenced by the choice of the horizon $M_n$, however, the selection rules MPVCcal and mPVCcal are the least affected.
We also investigated the influence  of the strength of dependence structure of design matrix $\X$ on the behaviour of  selection rules.
We studied the cases when the dependence between the covariates is respectively stronger and weaker than in the case described above. Namely  the covariances  $\Sigma_{X}(i,j)={0.8}^{|i-j|}$ and $\Sigma_{X}(i,j)=I\{i=j\}$ were considered.
For the above cases we took also different marginal variances of regressors equal to 0.5 and 2.
The error variance $\sigma^{2}$ was always set to one.
The experiments show that the probability of true model selection is smaller (and respective prediction error larger) than for initial scenario when the dependence is stronger or the variance of covariates larger. 
However, it turns out that the ranking of methods with respect to both considered measures remains the same in all above cases.
Experiments indicate also that for the considered selection criteria mean prediction error behaves approximately as a constant minus 
probability of a correct selection.\\
We also investigated the case of covariates $\bx^{n}_{l}$ having different distributions. Namely,  we considered the following regression scenario
\begin{equation*}
\Y=\betab'\Lb(\U)+\varepsilonb,
\end{equation*}  
where $\Lb(\cdot)=(L_1(\cdot),\ldots,L_{M_n}(\cdot))'$ is a vector consisting of  the consecutive orthonormal Legendre polynomials on $[-1,1]$ and $\U$ is random vector with continuous uniform distribution on $[-1,1]$. We considered the following list of models
\begin{enumerate}
\item[(L1)] $t=\{1,2,4\}$, $\betab=(1,1,1)'$
\end{enumerate}
with horizons $M_n=5,10,\ldots,25.$
The influence of the size of the list $M_n$ has been investigated. The sample size was set to $n=300$.
Figure 5 presents the results which are similar to that of the previous experiments indicating that mPVCcal and MPVCcal perform the
best in this case, and the second best is BIC.
\subsection{Real data example}
We consider {\tt bodyfat} data set  (\citet{Johnson1996}) consisting of records of  the percentage of fat in the body (dependent variable) together with  13 independent variables for $n=252$ individuals. Two independent variables were selected having the smallest p-values when the full linear model was fitted. They were abdomen and wrist circumference
and when used as predictors  resulted in the fitted model with a vector of estimated coefficients $\hat{\betab}=(0.7661,-2.8379)'$
and a variance of residuals $\hat{\sigma}^{2}=4.45$. A parametric bootstrap (see e.g. \citet{DavisonHinkley1997}) was employed to check how the considered selection criteria perform for this data set. Namely,  the true
model was the fitted linear model with the original two regressors, $\betab=\hat{\betab}$ and the normal errors with the variance equal to $\hat{\sigma}^{2}$. Additional superfluous explanatory variables were created
in pairs by drawing from the two-dimensional  normal distribution with independent components, which mean and  variance vector matched that of the original predictors. We considered $k=8,18,\ldots,58$ additional variables what amounted to horizons $M_n=10,20,\ldots,60$ when the true variables were accounted for.
 Thus $M_n/n$ ranged from $0.03$ to $0.23$.
 500 parametric bootstrap samples consisting of 252 observations each were created to mimic the original sample and 
 the considered selection criteria were employed to choose subset of potential $M_n$ variables. Figure 6 presents the results.
The results are similar to that of simulation experiments indicating that mPVCcal and MPVCcal perform the
best in this case, and the second best is BIC.
\begin{table}[h!]
\begin{scriptsize}
\begin{center}
\caption{Estimated probability of correct ordering based on $N=500$ trials.}
\begin{tabular}{llllllll}
\\\hline
Model &  $n=75$  & $n=100$ & $n=200$ & $n=300$ & $n=500$ & $n=1000$   \\ 
\hline
(M1) &0.16    &0.18  &0.39 &0.61 &  0.85 & 0.98\\
(M2) &0.69    &0.74   &0.91 &0.99 & 0.99 & 1\\
(M3) &0.99    & 1 & 1&1 & 1 & 1\\
(M4) &0.99    &  1&1 &1 & 1 & 1\\
\hline
Est. max. standard error $\leq$ 0.01&&&&\\
\hline
\label{TABORDERING}
\end{tabular}
\end{center}
\end{scriptsize}
\end{table}  
\newpage
\begin{figure}[h!]
\begin{center}$
\begin{array}{cc}
\includegraphics[scale=0.4]{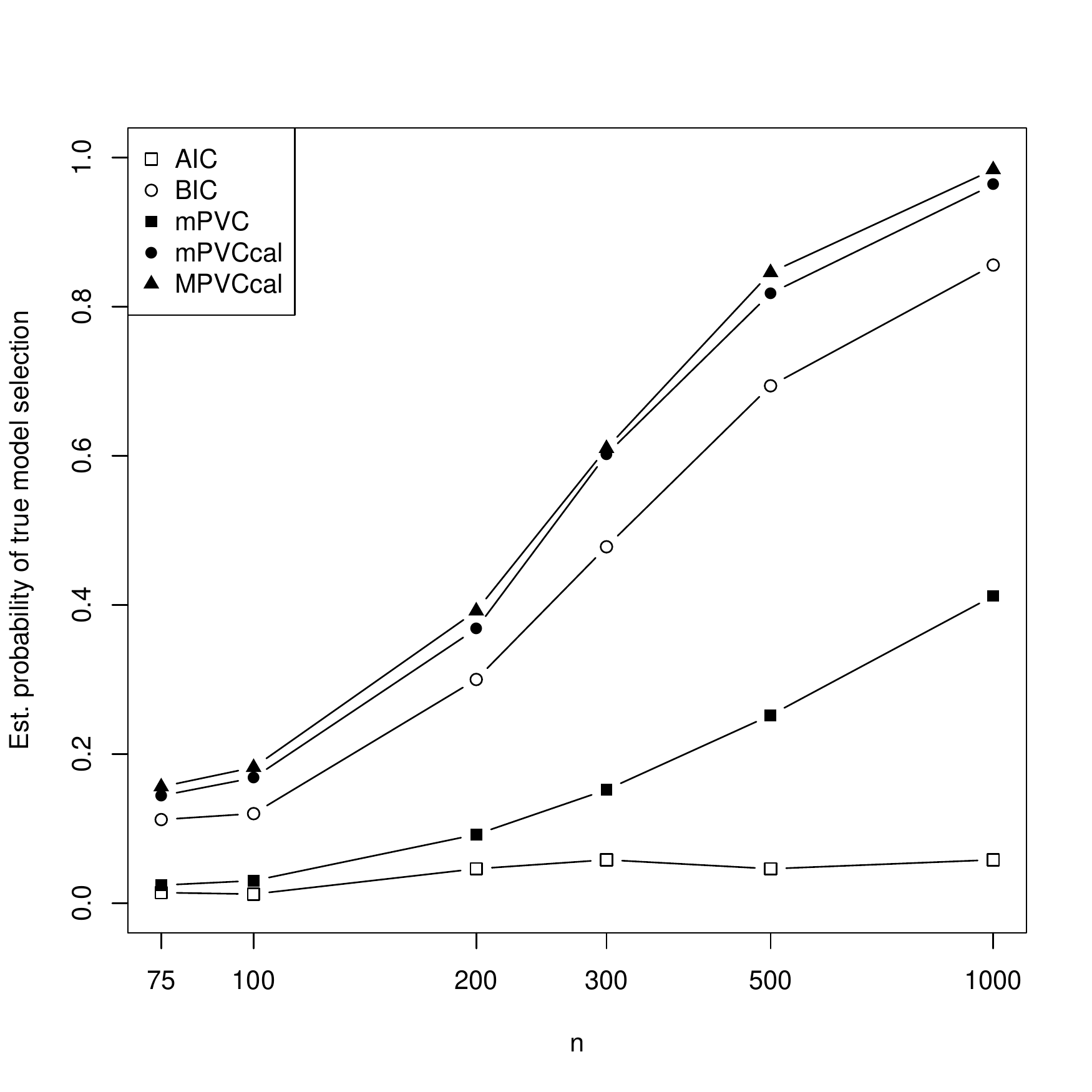} &
\includegraphics[scale=0.4]{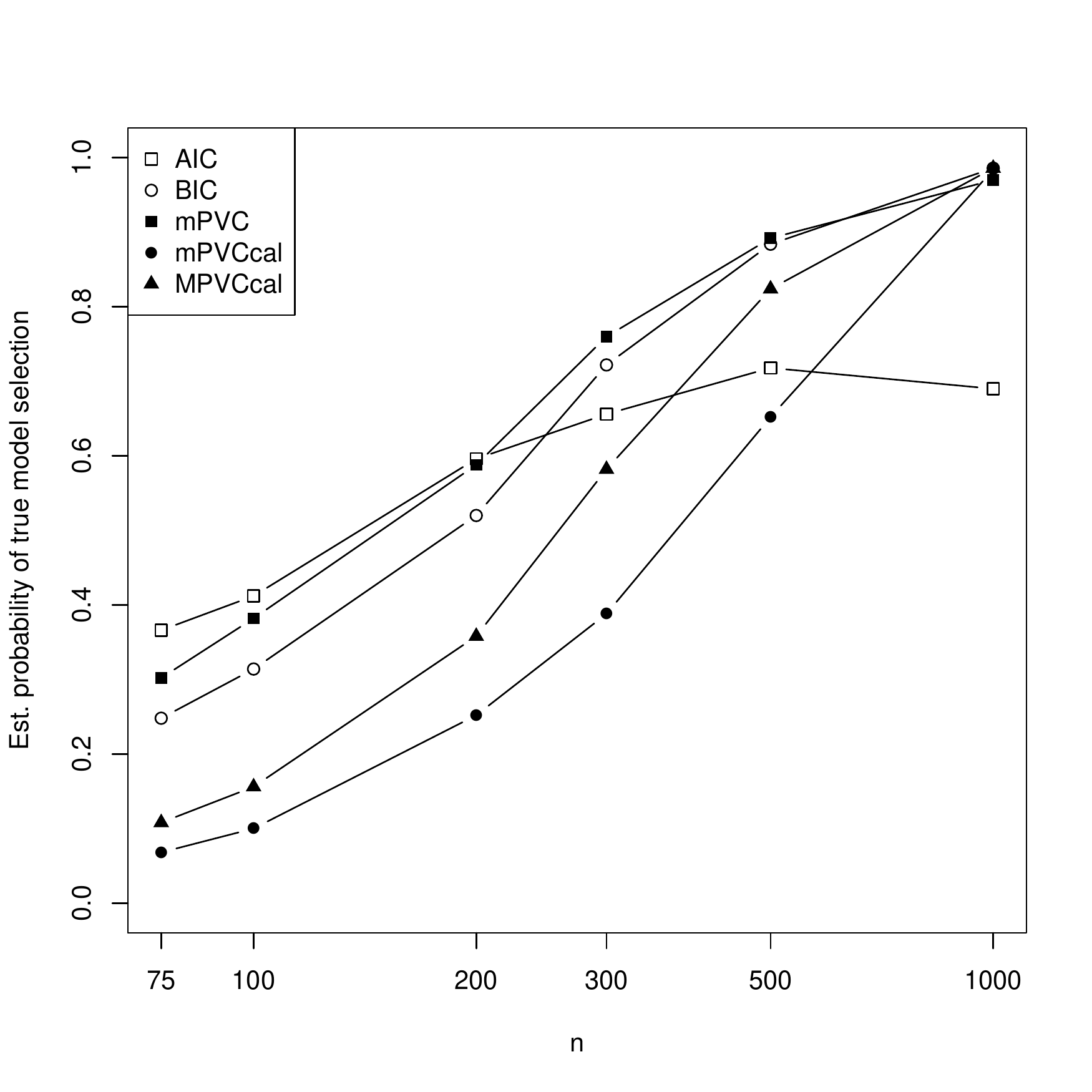} \\
(a)&(b)
\end{array}$
\end{center}
\end{figure}
\begin{figure}[h!]
\begin{center}$
\begin{array}{cc}
\includegraphics[scale=0.4]{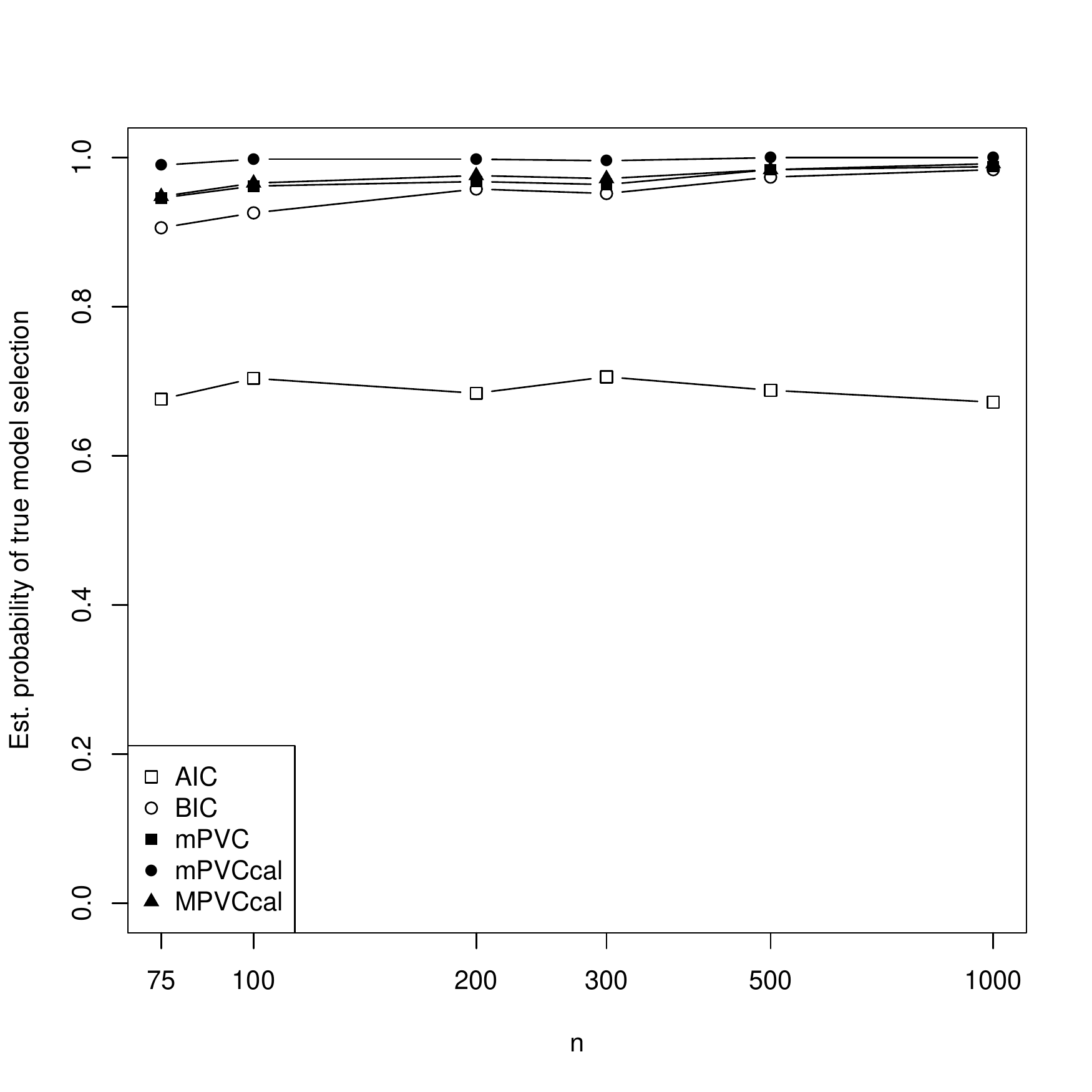} &
\includegraphics[scale=0.4]{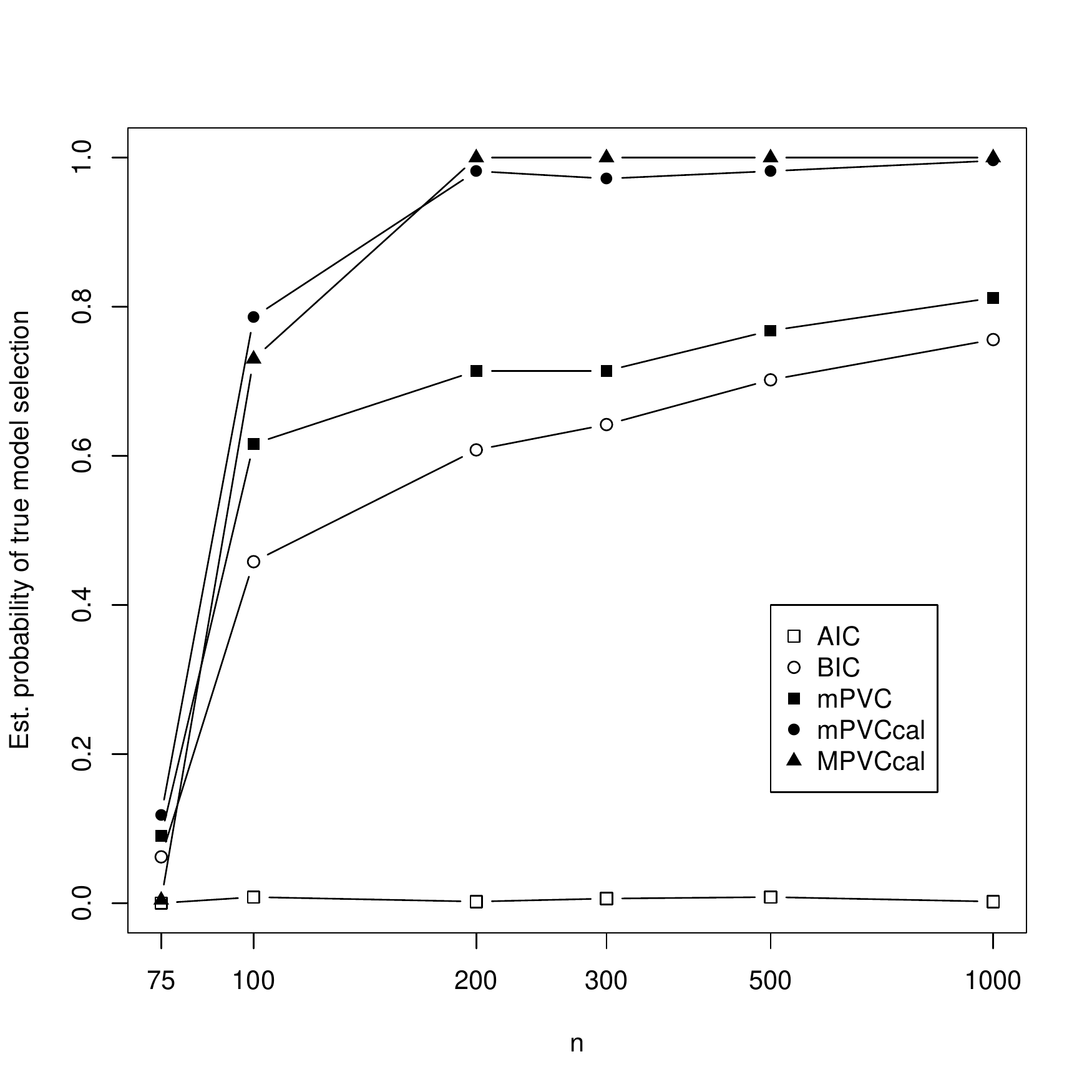} \\
(c)&(d)
\end{array}$
\end{center}
\caption{Estimated probabilities of correct model selection for models M1 (a), M2 (b), M3 (c) and M4 (d) with respect to $n$ (on a logarithmic scale) based on $N=500$ trials.}
\end{figure}
\newpage
\begin{figure}[h!]
\begin{center}$
\begin{array}{cc}
\includegraphics[scale=0.4]{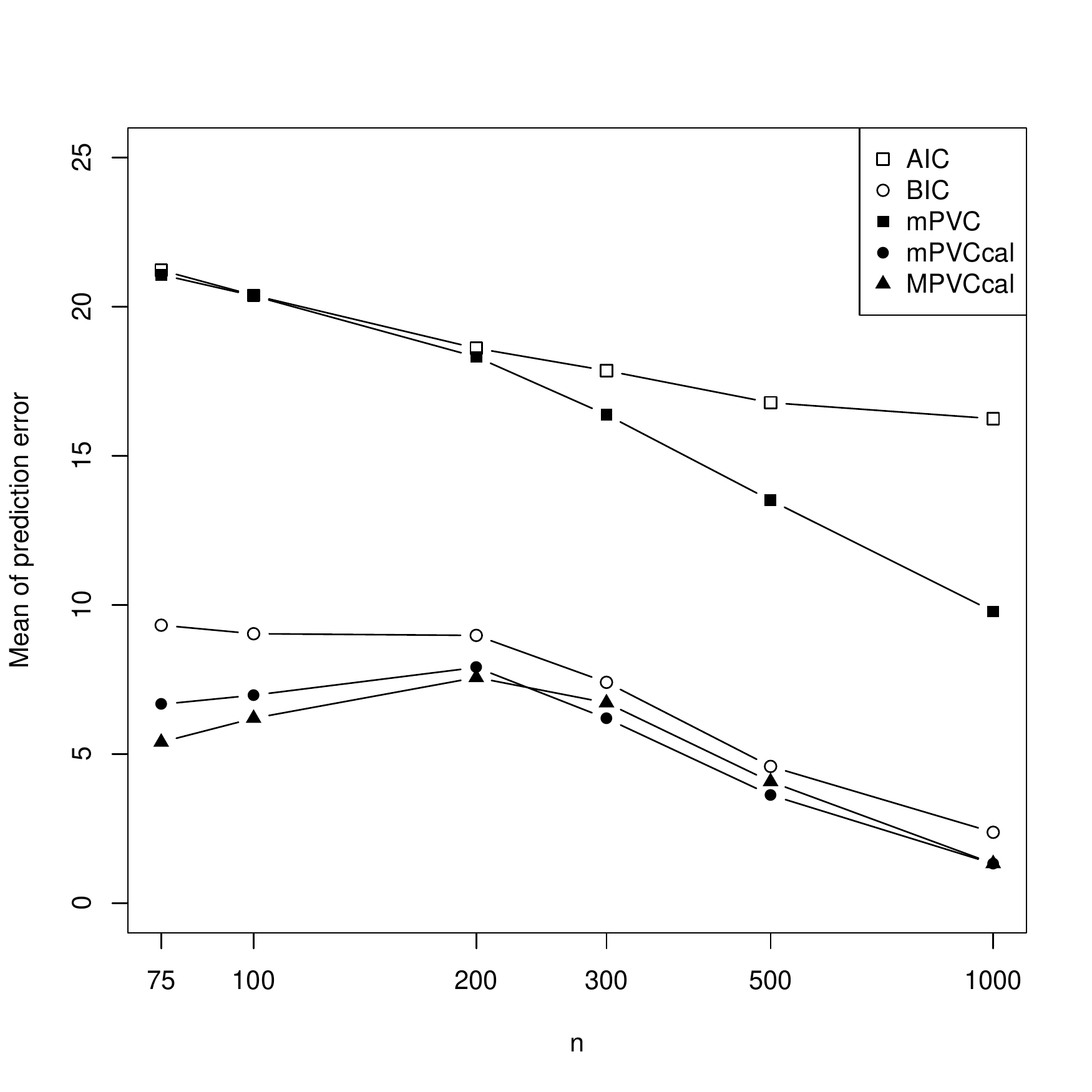} &
\includegraphics[scale=0.4]{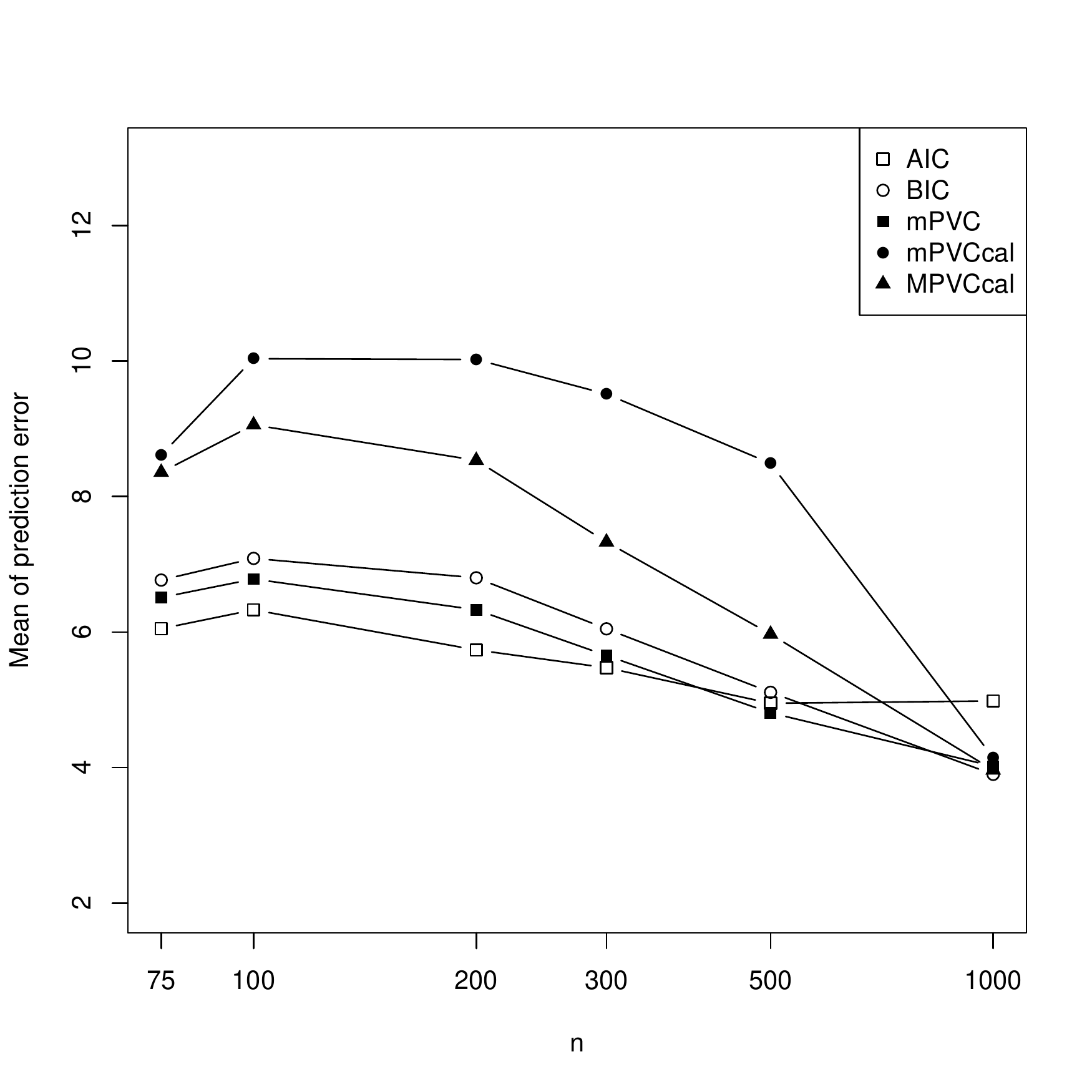} \\
(a)&(b)
\end{array}$
\end{center}
\end{figure}
\begin{figure}[h!]
\begin{center}$
\begin{array}{cc}
\includegraphics[scale=0.4]{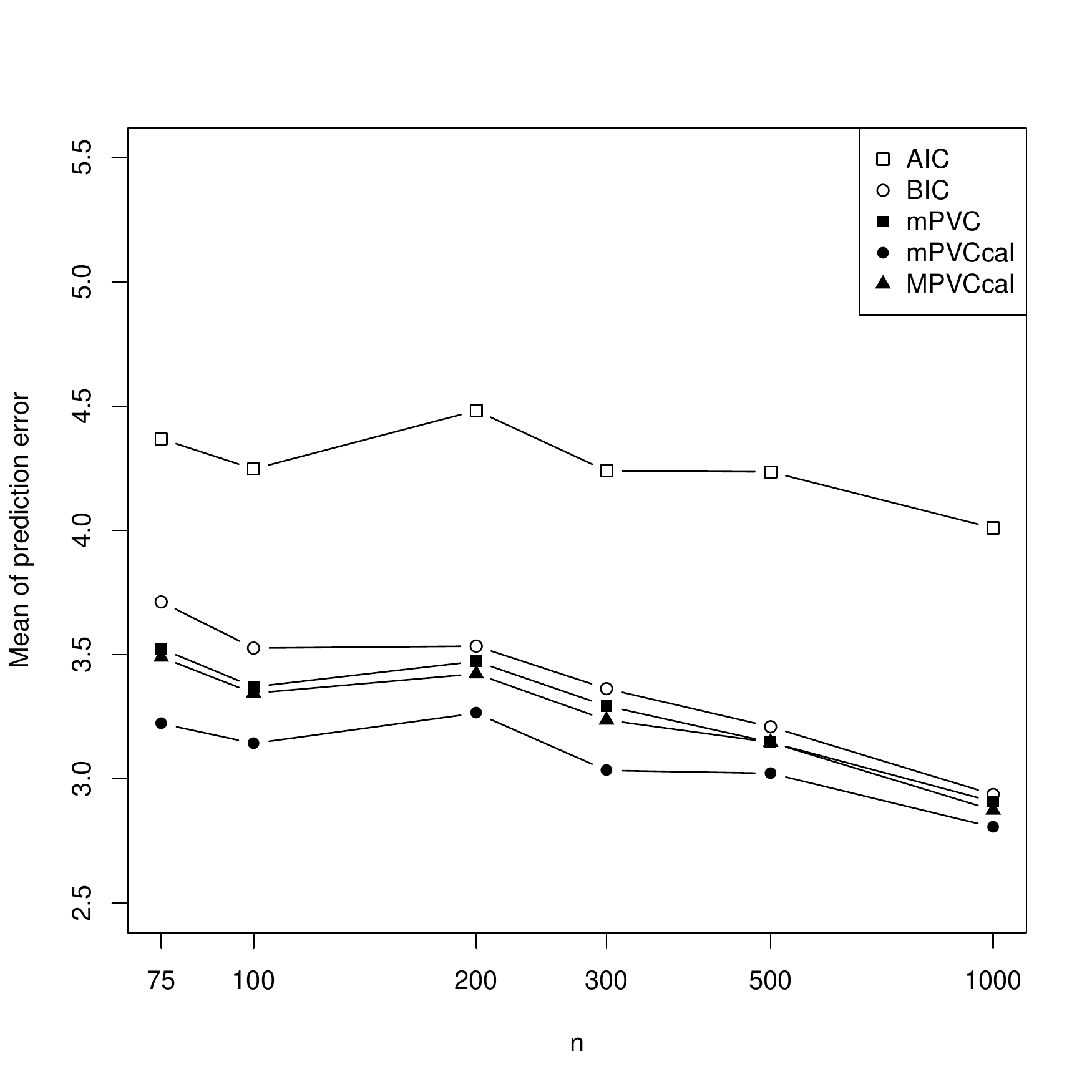} &
\includegraphics[scale=0.4]{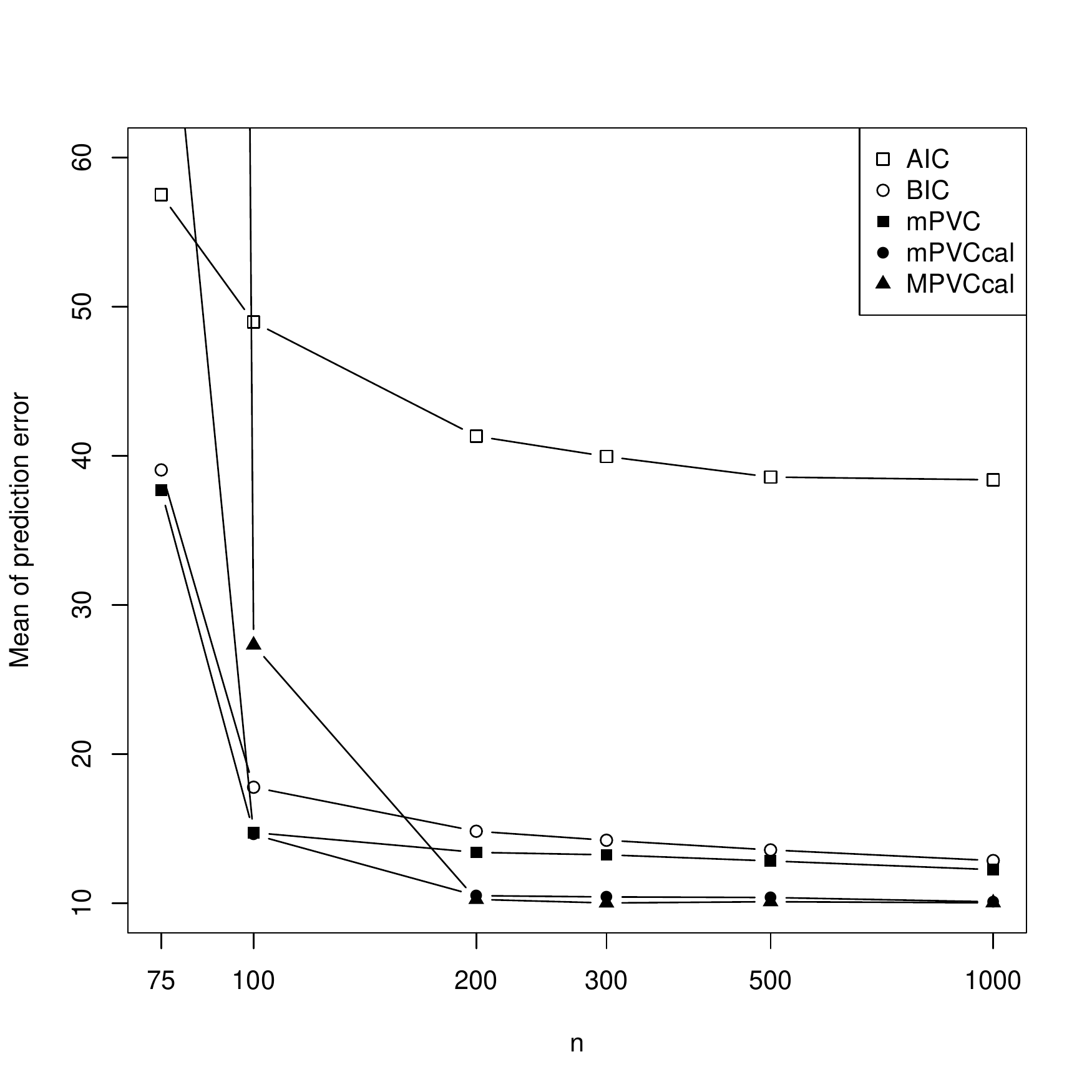} \\
(c)&(d)
\end{array}$
\end{center}
\caption{Means od prediction error for models M1 (a), M2 (b), M3 (c) and M4 (d) with respect to $n$ (on a logarithmic scale) based on $N=500$ trials.}
\end{figure}
\newpage
\begin{figure}[h!]
\begin{center}$
\begin{array}{cc}
\includegraphics[scale=0.4]{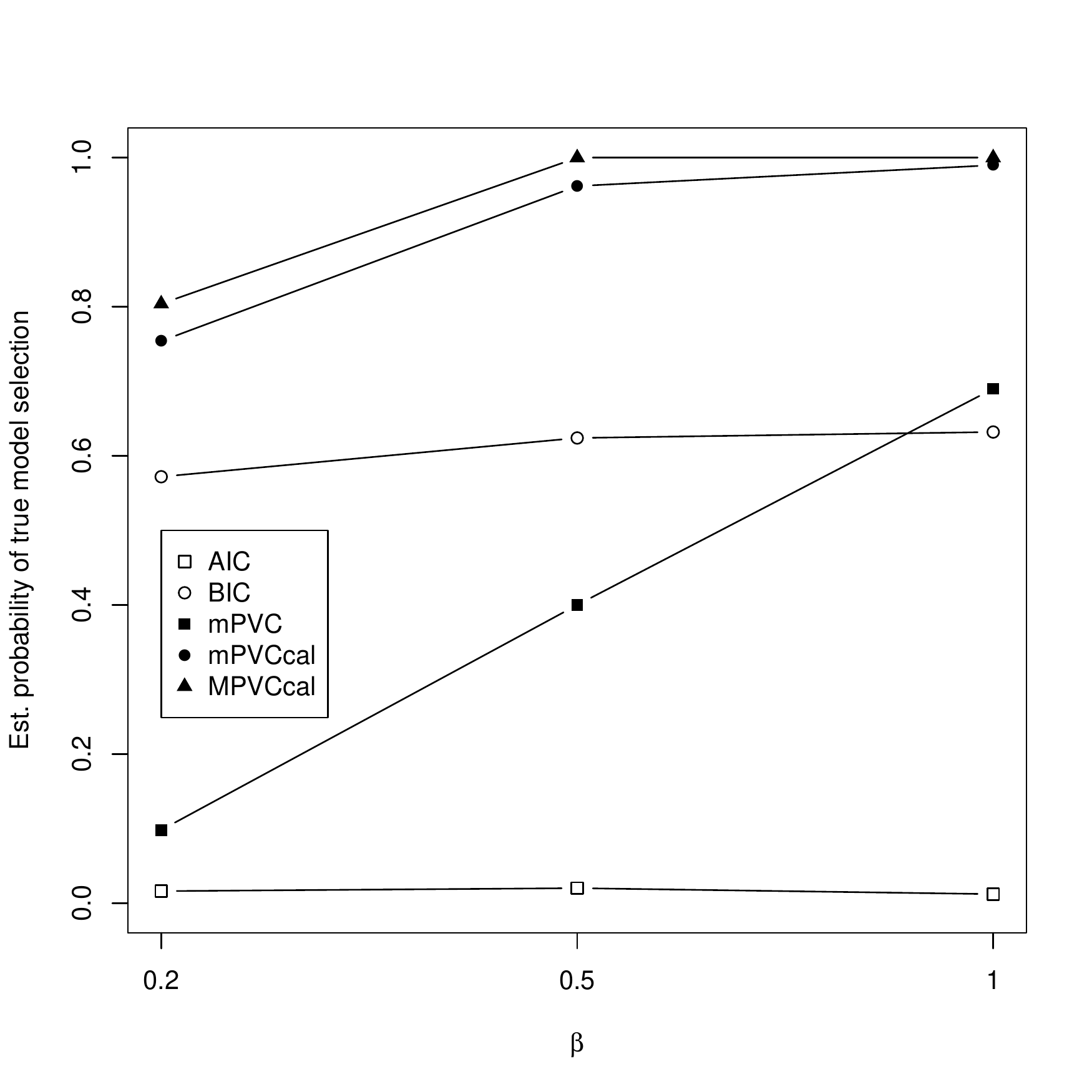} &
\includegraphics[scale=0.4]{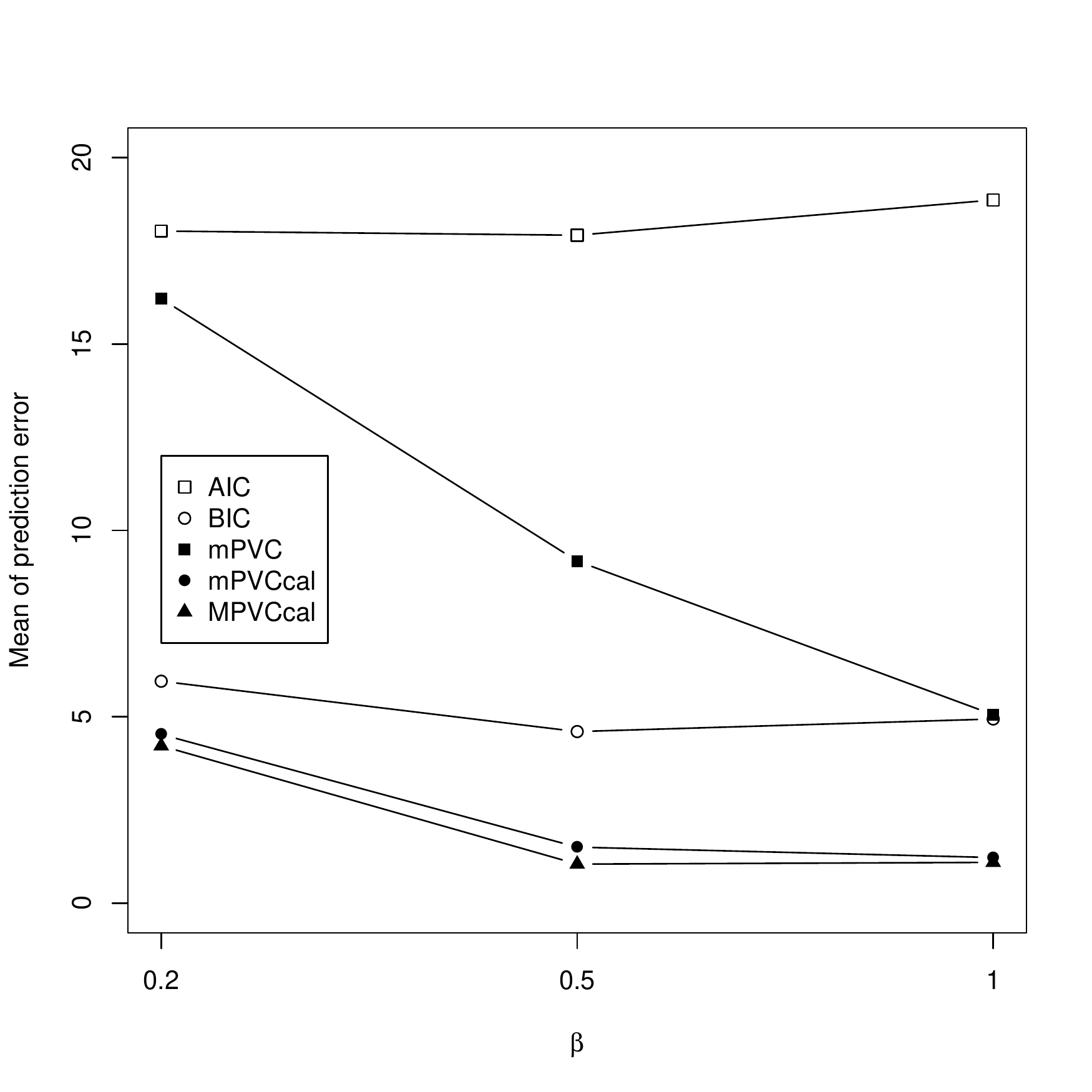} \\
(a)&(b)
\end{array}$
\end{center}
\caption{Estimated probabilities of correct model selection (a) and means of prediction error (b) with respect to value of parameter $\beta$ for model M1 for sample size $n=300$ based on $N=500$ trials.}
\end{figure}
\begin{figure}[ht!]
\begin{center}$
\begin{array}{cc}
\includegraphics[scale=0.4]{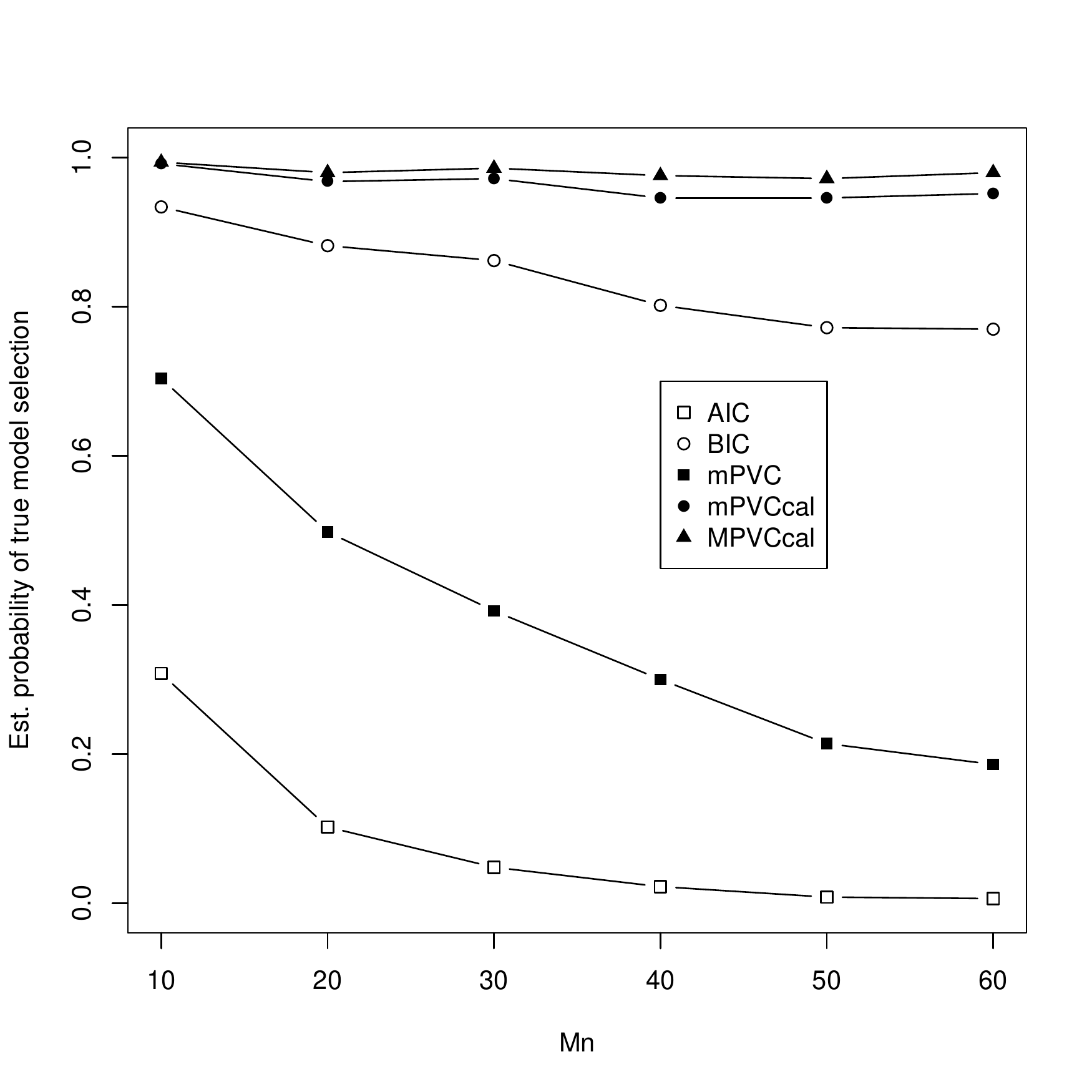} &
\includegraphics[scale=0.4]{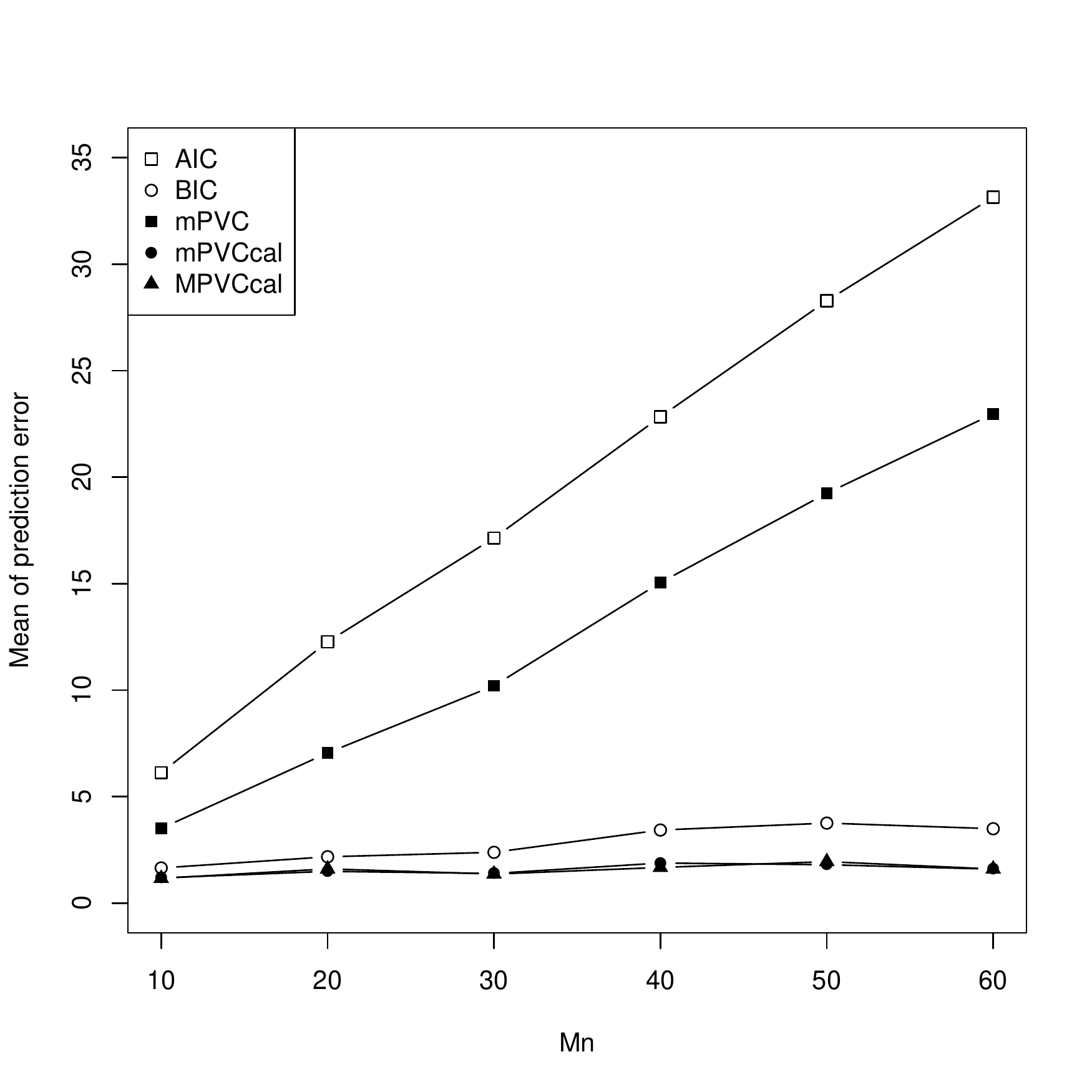} \\
(a)&(b)
\end{array}$
\end{center}
\caption{Estimated probabilities of correct model selection (a) and means of prediction error (b) with respect to $M_n$ for model M1 for sample size $n=1000$ based on $N=500$ trials.}
\end{figure}
\begin{figure}[ht!]
\begin{center}$
\begin{array}{cc}
\includegraphics[scale=0.4]{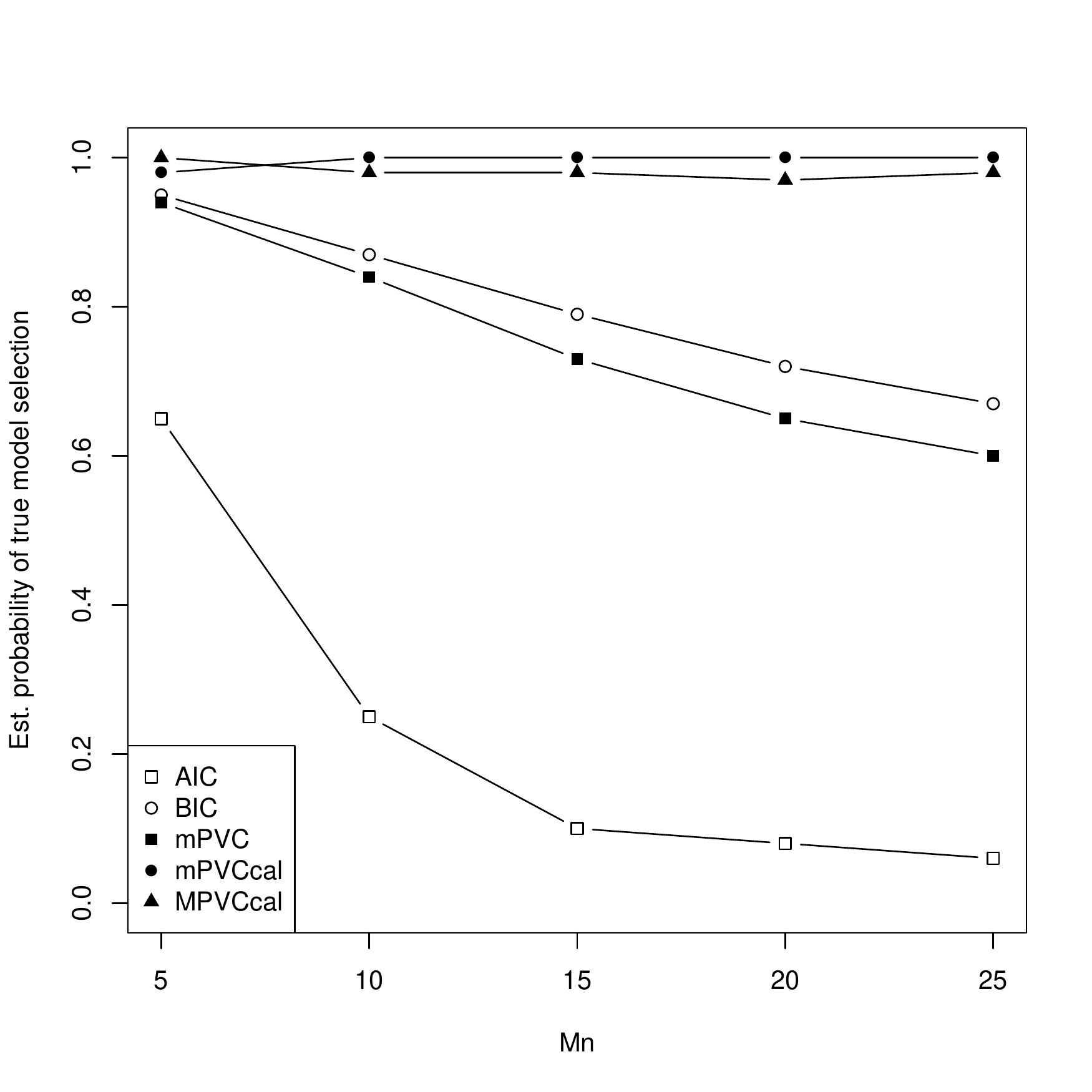} 
&\includegraphics[scale=0.4]{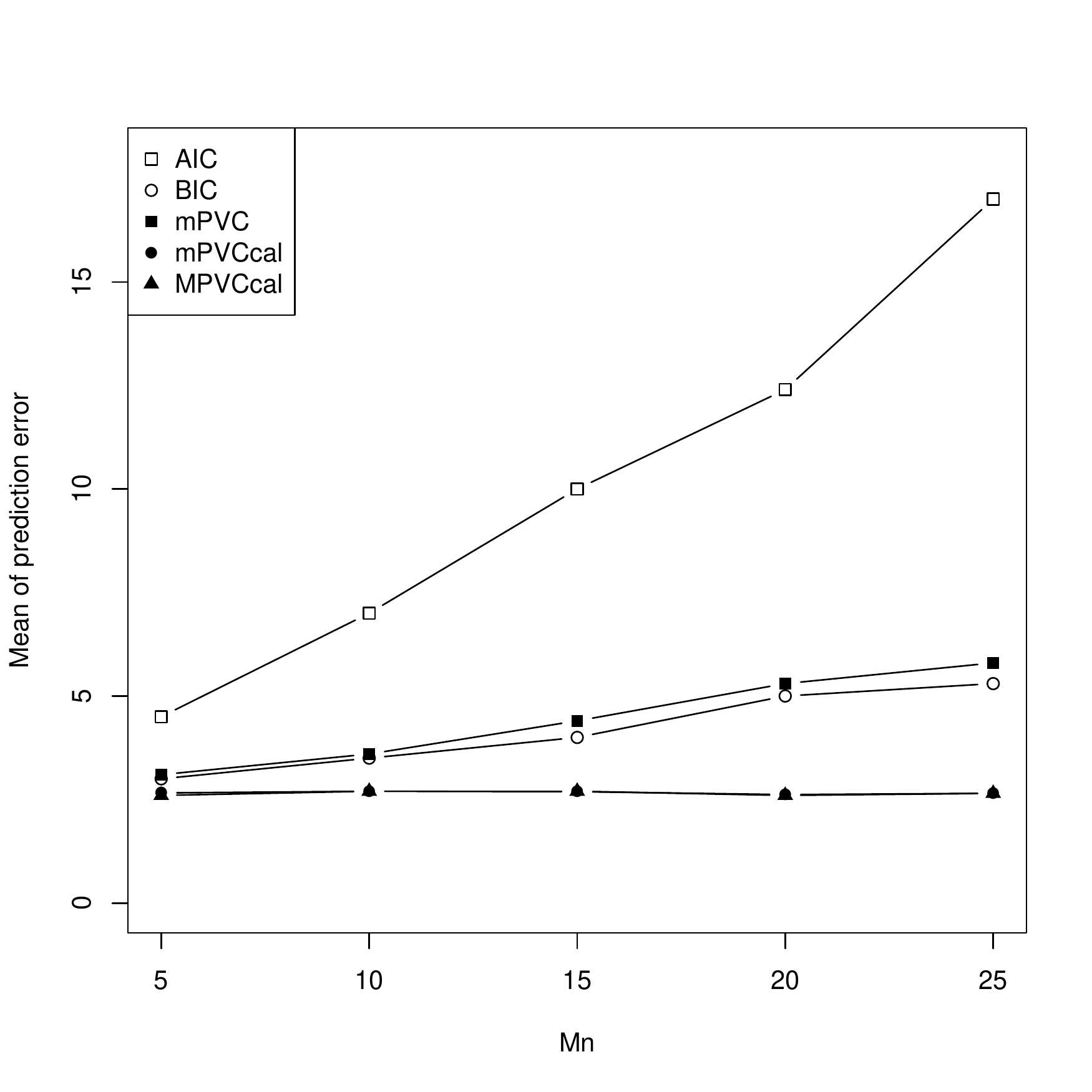} \\
(a)&(b)
\end{array}$
\end{center}
\caption{Estimated probabilities of correct model selection (a) and means of prediction error (b) with respect to $M_n$ for model (L1) based on $N=500$ trials.}
\end{figure}
\begin{figure}[ht!]
\begin{center}$
\begin{array}{cc}
\includegraphics[scale=0.4]{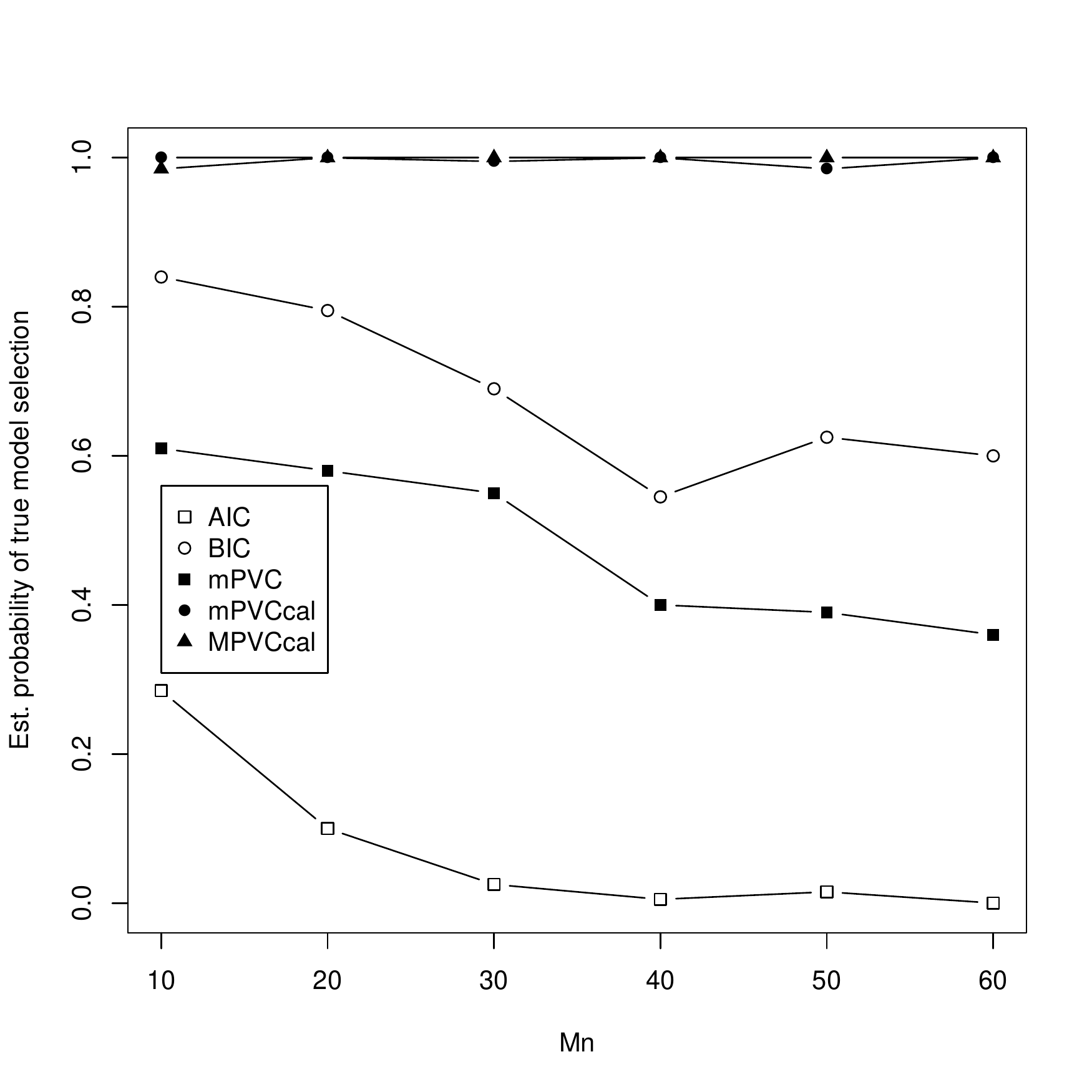} &
\includegraphics[scale=0.4]{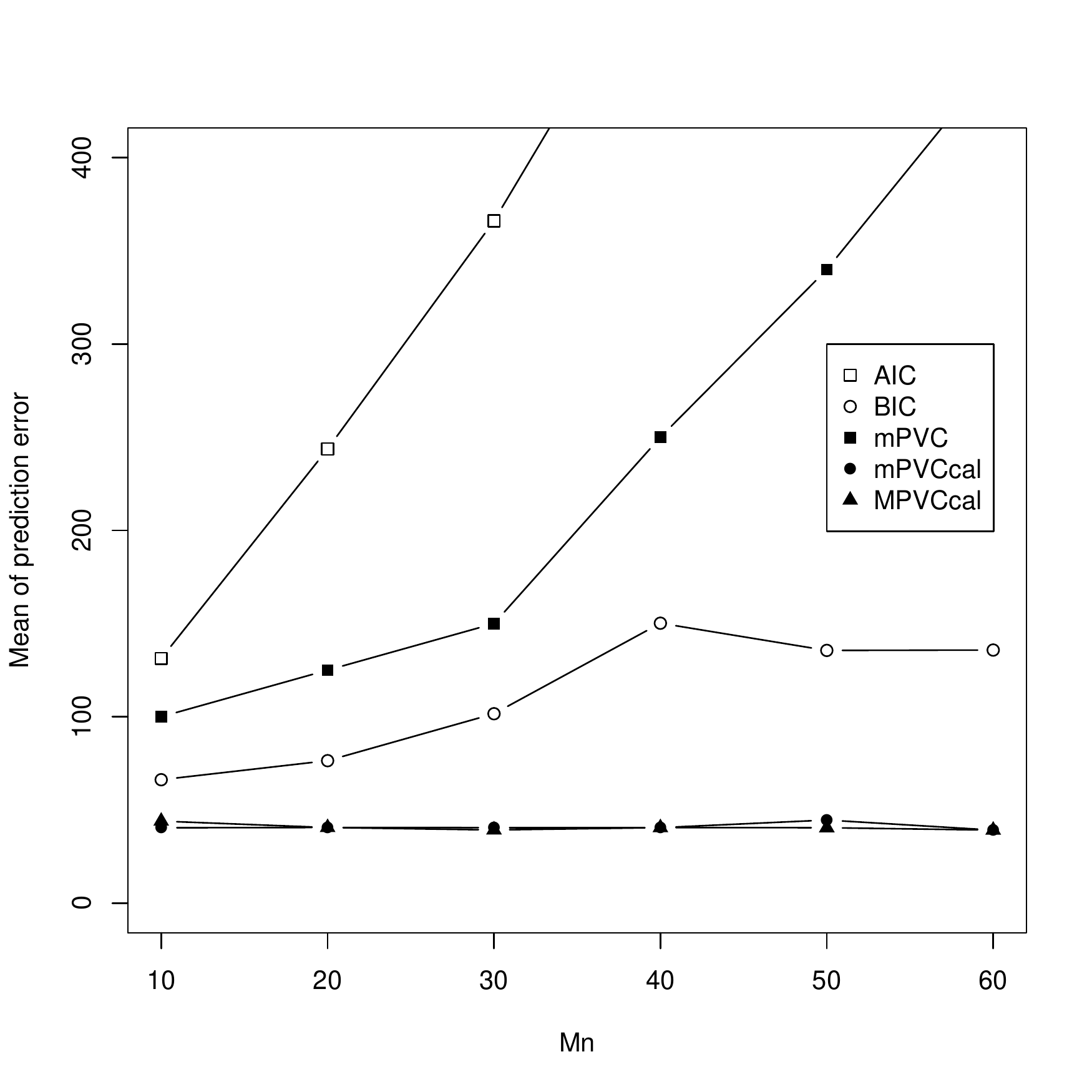} \\
(a)&(b)
\end{array}$
\end{center}
\caption{Estimated probabilities of correct model selection (a) and means of prediction error (b) with respect to $M_n$ for {\tt bodyfat} data set.}
\end{figure}
\newpage
\section{Appendix}
{\bf{ Proof of Lemma \ref{LemmaBeta} }}\\
The lemma is proved in \citet{PokarowskiMielniczuk2010}. For completeness we give an outline of proof here.
Recall that $B_{a,b}$ and $B(x,y)$ denote a random variable having  beta distribution with shape parameters $a$ and $b$ and beta function, respectively.
Let $B_{x}(a,b)=\int_{0}^{x}t^{a-1}(1-t)^{b-1}dt$ be the incomplete beta function.
It can be easily proved that
\begin{equation}
\label{LemmaBetaEq1}
aB_{x}(a,b)=x^{a}(1-x)^{b}+(a+b)B_{x}(a+1,b),
\end{equation}
and 
\begin{equation}
\label{LemmaBetaEq2}
B_{1-x}(b,a)=B(a,b)-B_{x}(a,b).
\end{equation}
Consider the case $a\geq 1$. Using (\ref{LemmaBetaEq1}), (\ref{LemmaBetaEq2}) and assumption $x>\frac{a-1}{a+b}$ we obtain the upper bound in (\ref{LemmaBetaeq1})
\begin{eqnarray*}
&&
P[B_{a,b}>x]=1-\frac{B_{x}(a,b)}{B(a,b)}=\frac{B_{1-x}(b,a)}{B(a,b)}=
\cr
&&
\frac{1}{B(a,b)b}\cdot(1-x)^{b}x^{a}[1+\frac{a+b}{b+1}(1-x)+\frac{(a+b)(a+b+1)}{(b+1)(b+2)}(1-x)^{2}+\ldots]\leq
\cr
&&
\frac{1}{B(a,b)b}\cdot(1-x)^{b}x^{a}[1+\frac{a+b}{b+1}(1-x)+\left(\frac{a+b}{b+1}\right)^{2}(1-x)^{2}+\ldots]=
\frac{(1-x)^{b}x^{a-1}}{B(a,b)b}(1+L(a,b,x)).
\end{eqnarray*}
In order to obtain the lower bound in (\ref{LemmaBetaeq1}) note that for $a\geq 1$
\begin{eqnarray*}
&&
\frac{1}{B(a,b)b}\cdot(1-x)^{b}x^{a}[1+\frac{a+b}{b+1}(1-x)+\frac{(a+b)(a+b+1)}{(b+1)(b+2)}(1-x)^{2}+\ldots]\geq
\frac{1}{B(a,b)b}\cdot(1-x)^{b}x^{a-1}.
\end{eqnarray*}
The case $a<1$ can be treated analogously.
\\
For ease of notation we assume in the following proofs that $\sigma^2=1$. Let $\Q(j)$ denote projection on the column space spanned by the regressors corresponding to coefficients in a given model $j$.\\
{\bf{ Proof of Lemma \ref{Auxilliary} }}\\
Consider first the case $j\subset t$.
Denote $\W=\E(\bx_{lt}\bx_{lt}')$, which in view of  assumption (A0) is positive definite.
Define $\Lambda_{n,j}=n^{-1}(\X\betab)'[\I-\Q(j)](\X\betab)>0$.
Let $\D_{j}$ be a $M_n\times j$ matrix of zeros and ones such that $\X\D_{j}$ consists of only these $j$ columns of $\X$ 
which correspond to model $j$. By assumption (A0) and  using the fact that $\X\betab=(\X \D_{t})\bar{\betab}$ where $\bar{\betab}=(\beta_{t_1},\ldots,\beta_{t_{p_t}})'$ we have $\Lambda_{n,j}\cp\lambda>0$ as $n\to\infty$.
The assertion follows from the fact that for $j\subset t$
\begin{equation}
\label{CONV}
n^{-1}(\X\betab)'[\I-\Q(j)](\X\betab)=
n^{-1}\bar{\betab}'\A\bar{\betab},
\end{equation} 
where
\begin{equation*} \A=[(\X\D_{t})'(\X\D_{t})]-[(\X\D_{t})'(\X\D_{t})]\bar{\D}_{j}[\bar{\D}_{j}'(\X\D_{t})'(\X\D_{t})\bar{\D}_{j}]^{-1}\bar{\D}_{j}'[(\X\D_{t})'(\X\D_{t})] 
\end{equation*}
and $\bar{\D}_{j}$ is a $p_t\times p_j$ matrix such that $\X\D_{j}=(\X\D_{t})\bar{\D}_{j}$.
Matrix $\W$ as a positive definite matrix can be decomposed as $\W=\W^{1/2}\W^{1/2}$ where $\W^{1/2}=\U\Xi^{1/2}\U '$, $\U$ is an orthogonal matrix and $\Xi$ is a diagonal matrix with positive diagonal.
 The right hand side of (\ref{CONV}) converges in probability to
\begin{eqnarray*}
&&
\lambda=\bar{\betab}'[\W-\W\bar{\D}_j({\bar{\D}_j}'\W\bar{\D}_j)^{-1}{\bar{\D}_j}'\W]\bar{\betab}=
\cr
&&
(\W^{1/2}\bar{\betab})'[\I-\W^{1/2}\bar{\D}_{j}({\bar{\D}_{j}}'\W\bar{\D}_{j})^{-1}{\bar{\D}_{j}}'(\W^{1/2})']\W^{1/2}\bar{\betab}>0
\end{eqnarray*}
since the columns of $\W^{1/2}$ are linearly independent. 
We have the following decomposition for $j\subset t$ 
\begin{equation}
\label{DECOMP}  
n^{-1}RSS(j)=n^{-1}\varepsilonb '(\I-\Q(j))\varepsilonb+
n^{-1}2(\X\betab)'(\I-\Q(j))\varepsilonb+\Lambda_{n,j}. 
\end{equation}
The first summand converges in probability to $\sigma^{2}$. The last summand $\Lambda_{n,j}\cp\lambda>0$, as has been already shown. 
Provided that ${\bf X'X}$ is invertible, $n^{-1}2(\X\betab)'(\I-\Q(j))\varepsilonb$ given $\X$ has
$N(0,v_n)$ distribution, where $v_n=n^{-1}\Lambda_{n,j}\cp 0$. Thus $n^{-1}2(\X\betab)'(\I-\Q(j))\varepsilonb\cp 0$.
This completes the first part of the proof. 
For $j\supseteq t$ the second and the third term in (\ref{DECOMP}) are equal to zero. This yields the second part of the assertion.\\ 
{\bf{ Proof of Lemma \ref{LemmaSupmodel} }}\\
Define $b_n=n(\exp(R_n/n)-1)$. It is easily seen that $b_n\geq R_n$ thus $b_n$ satisfies the condition imposed on $R_n$. For $M_n=p_t$
the assertion is obvious, thus we assume that $M_n>p_t$

We have the following inequality
\begin{eqnarray*}
&&
P\left\{n\log\left[\frac{RSS(t)}{RSS(f)}\right]>R_n\right\}=
P\left\{\frac{RSS(t)}{RSS(f)}>\exp\left(\frac{R_n}{n}\right)\right\}=
\cr
&&
P\{\varepsilonb '[\Q(f)-\Q(t)]\varepsilonb > {b}_n n^{-1}\varepsilonb '[\I-\Q(f)]\varepsilonb\}\leq
\cr
&&
P\{\varepsilonb '[\Q(f)-\Q(t)]\varepsilonb > {b}_n n^{-1}(n-M_n-d_n)\} +
\cr
&& 
P\{\varepsilonb '[\I-\Q(f)]\varepsilonb \leq n-M_n-d_n\},
\end{eqnarray*}
where $d_n = (n-M_n)^{(1+\delta)/2}$, for some $\delta\in(0,1)$.
Matrix ${\bf X'X}$ has rank $M_n$ and it follows 
 that $\varepsilonb '[\Q(f)-\Q(t)]\varepsilonb\sim\chi^{2}_{M_n-p_t}$ and $\varepsilonb '[\I-\Q(f)]\varepsilonb\sim{\chi^{2}}_{n-M_n}$ 
(since $\sigma^2=1$). By an inequality for cumulative distribution function of a chi-square distribution,
\begin{equation*}
P(\chi^{2}_{k}\leq k - \delta_0)\leq \exp\{-(4k)^{-1}\delta_{0}^{2}\},
\end{equation*}
for $\delta_0 > 0$ (see \cite{Shibata1981}). Thus we have
\begin{equation*}
P\{\varepsilonb '[\I-\Q(f)]\varepsilonb\leq n-M_n-d_n\}\leq 
\exp\left[-\frac{d_{n}^2}{4(n-M_n)}\right]\to 0,
\end{equation*}
as $n\to\infty$, since $M_n/n\to 0$.
Let $\gamma_n=b_n(1-M_n/n-d_n/n)$. As $\varepsilonb '[\Q(f)-\Q(t)]\varepsilonb\sim \chi^2_{M_n-p_t}$ by Chebyschev
inequality we have  
\begin{equation*}
P\{\varepsilonb '[\Q(f)-\Q(t)]\varepsilonb -(M_n-p_t) > \gamma_n -(M_n-p_t) \}\leq \frac{2(M_n-p_t)}{[{\gamma}_n -(M_n-p_t)]^2}\to 0,
\end{equation*}
where the last convergence follows from $(\gamma_n-M_n)/\sqrt{M_n}\to \infty$. This completes the proof.\\
{\bf{ Proof of Lemma \ref{LemmaSubmodel} }}\\
In view of conditions (A1.3) and (A1.4) matrix $(\X '\X)^{-1}$ exists with probability tending to one (see the proof of Theorem 2 in \citet{ZhengLoh1997}).
Recall that $T_k$ is a t-statistic corresponding to the $k$th variable.
It suffices to prove that for any $c_n\to 0$ $P[\min_{i\in t}\log(RSS(f-\{i\})/RSS(f))<c_n]\to 0$.
Noting that 
\begin{eqnarray*}
&&
\frac{RSS(f-\{i\})}{RSS(f)}=\frac{T_{i}^{2}}{n-M_n}+1, 
\end{eqnarray*}
we obtain that 
\begin{eqnarray*}
&&
P[\min_{i\in t}\log\frac{RSS(f-\{i\})}{RSS(f)}< c_n]\leq
P[\min_{i\in t}T_{i}^{2}<(n-M_{n})(\exp(c_n)-1)]
\cr 
&&\leq P(\min_{i\in t}T_{i}^{2}<(n-M_n)(\exp(c_n)-1)).
\cr
&& 
\end{eqnarray*}
Since $\exp(c_n)-1=c_n+o(c_n)$ it suffices to show that $P[\min_{i\in t}T_{i}^{2}<Cnc_n]\to 0$, for some $C>0$. 
This follows from the proof of Theorem 2 in \citet{ZhengLoh1997} who  proved that under conditions of this Lemma 
$P[\min_{i\in t}\hat{\sigma}^{2}T_{i}^{2}<nc_n]\to 0$, for any $c_n$ such that  $c_n\to 0$.
 Now the  required convergence  follows from the fact that $\hat{\sigma}^{2}\cp\sigma^{2}$.
\bibliography{References}\bibliographystyle{plainnat}

\end{document}